\documentclass[11pt]{amsart}
\usepackage{amsmath,amsthm, amscd, amssymb, amsfonts}

\usepackage{hhline}

\newcommand{\fpd}{\operatorname{ FP-dim} \,}
\newcommand{\qdim}{\operatorname{\bf dim} \,}

\DeclareMathOperator*{\prode}{\boxtimes}

\newcommand{\ojo}{\ \\ \hspace*{-1cm}\textbf{ojo!}\hspace{2mm}}

\newcommand\boxe{\begin{tabular}{|p{0,1cm}|}
\hline \\ \hline \end{tabular}}

\newcommand\boxee{\begin{tabular}{|p{0,3cm}|}
\hline \\ \hline \end{tabular}}

\newcommand\boxu{\begin{tabular}{|p{0,1cm}|}
\\ \hline  \end{tabular}}

\newcommand\boxp{\begin{tabular}{|p{0,1cm}|}
 \hline \\ \end{tabular}}

\newcommand\boxa{\begin{tabular}{p{0,1cm}|}
\hline \\ \hline \end{tabular}}

\newcommand\boxc{\begin{tabular}{|p{0,1cm}}
\hline \\ \hline \end{tabular}}

\newcommand\mvert[2]{\begin{tiny}\begin{matrix}#1\vspace{-4pt}\\#2\end{matrix}\end{tiny}}
\DeclareMathOperator*{\Tim}{\times}
\newcommand{\Times}[2]{\sideset{_#1}{_#2}\Tim}

\newcommand{\cx}{{\daleth}}
\newcommand{\Y}{{\mathcal Y}}
\newcommand{\R}{{\mathcal R}}

\newcommand{\daga}{{\dagger}}

\newcommand{\ku}{\Bbbk}

\newcommand{\Z}{{\mathbb Z}}

\newcommand{\G}{{\mathcal G}}

\newcommand{\Q}{{\mathcal Q}}

\newcommand{\C}{{\mathcal C}}

\newcommand{\Ec}{{\mathbf E}}
\newcommand{\uno}{{\bf 1}}

\newcommand{\Dc}{{\mathbf D}}
\newcommand{\B}{{\mathcal B}}
\newcommand{\T}{{\mathcal T}}
\newcommand{\Hc}{{\mathcal H}}
\newcommand{\Vc}{{\mathcal V}}
\newcommand{\Pc}{{\mathcal P}}
\newcommand{\Ee}{{\mathcal E}}

\newcommand{\Ss}{{\mathcal S}}
\newcommand{\Vect}{\operatorname{Vec}}
\newcommand{\End}{\operatorname{End}}

\newcommand\tr{\operatorname{tr}}
\newcommand\Rep{\operatorname{Rep}}

\newcommand{\unosigma}{{\bf 1}^{\sigma}}

\newcommand{\fde}{{\triangleright}}
\newcommand{\fiz}{{\triangleleft}}

\newcommand{\prin}{t}
\newcommand{\fin}{b}
\newcommand{\pri}{r}
\newcommand{\fine}{l}

\newcommand\V{\operatorname{Vec}}

\numberwithin{equation}{section}\theoremstyle{plain}

\newtheorem{theorem}{Theorem}[section]
\newtheorem{lema}[theorem]{Lemma}
\newtheorem{cor}[theorem]{Corollary}

\newtheorem{prop}[theorem]{Proposition}

\theoremstyle{definition}
\newtheorem{definition}[theorem]{Definition}
\newtheorem{exa}[theorem]{Example}

\theoremstyle{remark}
\newtheorem{obs}[theorem]{Remark}
\newtheorem{rmk}[theorem]{Remarks}
\newtheorem{step}{Step}

\newcommand\id{\operatorname{id}}
\newcommand\idd{\mathbf{id}}
\newcommand\iddh{\mathbf{id}}
\newcommand\iddv{\mathbf{id}}
\newcommand\op{\operatorname{op}}

\def\pf{\begin{proof}}
\def\epf{\end{proof}}

\theoremstyle{remark}

\begin{document}

\renewcommand{\baselinestretch}{1.2}

\thispagestyle{empty}
\title{Tensor categories attached to double groupoids}
\author[Andruskiewitsch and Natale]{ Nicol\'as Andruskiewitsch and Sonia Natale}
\address{\noindent
Facultad de Matem\'atica, Astronom\'\i a y F\'\i sica,
Universidad Nacional de C\'ordoba.  CIEM -- CONICET.
(5000) Ciudad Universitaria, C\'ordoba, Argentina}
\email{andrus@mate.uncor.edu, \newline \indent \emph{URL:}\/
http://www.mate.uncor.edu/andrus} \email{natale@mate.uncor.edu,
\newline \indent \emph{URL:}\/ http://www.mate.uncor.edu/natale}
\thanks{This work was partially supported by CONICET, Fundaci\' on
Antorchas, Agencia C\'ordoba Ciencia, ANPCyT    and Secyt (UNC)}
\date{\today}
\begin{abstract} The construction of a quantum groupoid out of
a double groupoid satisfying a filling condition and a
perturbation datum is given. Several important classes of examples
of tensor categories are shown to fit into this construction.
Certain invariants such as a pivotal group-like element and
quantum and Frobenius-Perron dimensions of simple objects are
computed.
\end{abstract}

\maketitle

\setcounter{tocdepth}{1} \tableofcontents

\section*{Introduction}

 The main goal of this paper is the construction of a large
class of examples of weak Hopf algebras, also called \emph{quantum
groupoids}, from a  certain, quite general, class of double
groupoids.

\medbreak Quantum groupoids have been introduced not long ago \cite{bsz, bnsz}.
They are a natural generalization of the notion of a groupoid in a
non-commutative context. Heuristically, a finite quantum
groupoid consists of two  algebra structures in a finite
dimensional vector space, subject to a set of compatibility
conditions, and possessing an antipode.

 Quantum groupoids are interesting objects due to the fact
that they give rise, through its representation theory, to rigid
tensor categories. Tensor categories are important in several
areas of mathematics and physics. In particular, semisimple finite
quantum groupoids give rise to
semisimple rigid tensor categories with some finiteness
conditions. The key fact about quantum groupoids is that this
property does not only make them into a tool for constructing
tensor categories, but \emph{every} fusion category is the
representation category of a quantum groupoid, in view of results
of Hayashi and Ostrik.

\medbreak A double groupoid is a groupoid object in the category of
groupoids. Double groupoids were originally introduced by
Ehresmann \cite{ehr} in the early sixties, and later studied by a number of mathematicians interested in the search of a non commutative
relative higher homotopy groupoid of a topological space. Several
results on double groupoids have appeared in the literature since
then. Remarkably, their connection with crossed modules and a
higher Van Kampen Theory has been established in the work of
Brown, Higgins, Spencer,  {\it et al}. See for instance the survey
paper \cite{brown} and references therein.

A double groupoid can be roughly understood as a set of 'boxes'
with two groupoid compositions -the \emph{vertical} and
\emph{horizontal} compositions-, together with coherent groupoid
compositions of the sides, such that the boxes compositions obey a
set of compatibility conditions, one  of them being the so called
\emph{interchange law}.

\medbreak It seems natural to consider the following construction:
given a (finite) double groupoid $\T$, take the vertical and
horizontal groupoid algebra structures on the vector space spanned
by the boxes of $\T$. This construction was considered by the
authors in \cite{AN}; the necessary and sufficient condition for
this to produce a quantum groupoid is the \emph{vacancy} of $\T$.
The results in \cite{AN} gave a generalization of a celebrated
construction in Hopf algebra theory, studied by several people,
including G. I. Kac, Takeuchi and Majid. Essentially, a
\emph{vacant} double groupoid corresponds to an \emph{exact
factorization} of a groupoid. The resulting quantum groupoid is in
this case a kind of abelian `bicrossed product'.

\medbreak The main result of this paper is the determination of a
quantum groupoid structure in the span of the boxes in a double
groupoid satisfying a natural \emph{filling condition}, by introducing a certain perturbation in (any) one of
the groupoid algebra structures. This is achieved by considering a
family of 'corner' functions defined on a double groupoid. In fact,
the perturbation is done for a more general class of
functions. See Theorems \ref{bicross}, \ref{cx}.

Unlike in the vacant context, this more general construction does
\emph{not} fit into any known bicrossed product construction coming from matched pairs. It
is apparent that our approach is related to the formalism of
Ocneanu cells; however, the eventual precise relation still
remains to be developed.

We show that several important classes of examples of tensor
categories fit into our construction. For instance, for every
separable algebra $R$, the tensor category of $R$-bimodules (with
tensor product $\otimes_R$) is isomorphic to the category of
representation of a (not canonical) quantum groupoid arising from
a double groupoid.

\medbreak The source and target subalgebras of the resulting
quantum groupoid turn out to be isomorphic to the groupoid algebra
of the \emph{core groupoid} of $\T$. The core groupoid
of a double groupoid has been studied by Brown and Mackenzie: it
is known that it determines the whole double groupoid, under
appropriate restrictions. See \cite[Theorem 2.7]{BM}. However, as pointed out by Brown and Mackenzie in {\it loc. cit.}, it is not expected that general double groupoids can be described in terms of 'more familiar' structures.

\medbreak The intrinsic combinatorics of double groupoids reflect into 'Hopf-theo\-retic' features of the associated quantum groupoids. An instance of this principle is illustrated in Example \ref{no-siempre}, where we show that for any positive rational number $r$, there exists a double groupoid for which the square of the antipode of the associated weak Hopf algebra has $r$ as one of its eigenvalues.

\medbreak The paper is organized as follows. Section \ref{do-gpd}
contains the basic facts about double groupoids needed for the
construction. It introduces the corner functions and studies its
main properties. Sections \ref{wha} and \ref{extcocycles} present
the main results of the paper: the construction of the quantum
groupoid from a double groupoid, and a deformation of  this
construction via a certain cohomological data. Several properties
of these quantum groupoids are also studied in these sections.
Finally we consider in Section \ref{ejemplos} various examples of
tensor categories which fit into our construction.

\bigbreak {\bf Conventions.} Along this paper, in the case where
$f, g$ are composable arrows in a groupoid, their composition will
be indicated by juxtaposition from left to right, that is, we
shall use the notation $fg$ instead of $gf$.

For a groupoid $\G \rightrightarrows \Pc$, with base $\Pc$, we shall
identify $\Pc$ with a subset of $\G$ via the identity map $\G \to \Pc$;
when there is no ambiguity we shall speak of 'the groupoid $\G$' instead of
$\G \rightrightarrows \Pc$.

\medbreak {\bf Thanks.} A preliminary stage of this research  was
conducted at the University of Miskatonic (Arkham); both authors
are grateful to C. D. Ward and H. West for the kind hospitality.

Results of this paper were reported by the first-named author at
the Ferrara Algebra Workshop, June 16th-19th 2004. He thanks
Claudia Menini for the generous invitation.  They were also
reported by the second-named author at the International
Conference on Quantum Groups, Haifa, 05/07--12/07 2004; she thanks
the organizers for the kind invitation and hospitality. Both
authors thank the referee for his/her interesting comments.

\section{Double groupoids}\label{do-gpd}

\subsection{Definition of double groupoid}\label{recalling}

\

A (finite) \emph{double groupoid}  $\T$ is a  groupoid
object in the category of (finite) groupoids. It is customary to
represent a double groupoid in the form of four related groupoids
$$\begin{matrix} \B &\rightrightarrows &\Hc
\\\downdownarrows &&\downdownarrows \\ \Vc &\rightrightarrows &\Pc \end{matrix}$$
subject to a set of axioms. See \cite{ehr, bs}. Throughout this
paper, we shall keep the conventions and notations from
\cite[Section 2]{AN}.

The source and target maps  of these groupoids are indicated by
$\prin, \fin: \B \to \Hc$; \, $\pri, \fine: \B \to \Vc$; \,
$\pri, \fine: \Hc \to \Pc$; \, $\prin, \fin: \Vc \to \Pc$
(`top', `bottom', `right' and `left'). An element $A\in \B$ is
depicted as a box $$A =
\begin{matrix} \quad t \quad \\ l \,\, \boxe \,\, r \\ \quad b\quad
\end{matrix}$$ where $t(A) = t$, $b(A) = b$, $r(A) = r$, $l(A) = l$,
and the four vertices of the square representing $A$ are $tl(A) = lt(A)$,
$tr(A) = rt(A)$, $bl(A) = lb(A)$, $br(A) = rb(A)$. A box $A\in \B$
is, in general, \emph{not} determined by its boundary.

According to our conventions, horizontal and vertical composition
of boxes will be written from left to right and from top to
bottom, respectively. We shall write $A \vert B$ if $r(A)=l(B)$
(so that $A$ and $B$ are horizontally composable), and
$\displaystyle\frac{A}{B}$ if $b(A)=t(B)$ (so that $A$ and $B$ are
vertically composable). The notation $AB$ (respectively,
$\begin{matrix}A \vspace{-4pt}\\B\end{matrix}$) will indicate the
horizontal (respectively, vertical) composition; this notation
will always implicitly assume that $A$ and $B$ are composable in
the appropriate sense. Sometimes we will denote a box having an
identity in the top as \begin{tabular}{|p{0,1cm}|} \hhline{|=|} \\
\hline\end{tabular}\,; or, if the box has an identity in the left,
as \begin{tabular}{||p{0,1cm}|} \hline \\
\hline \end{tabular}\,, etc.

Compositions verify the following. Let $A = \begin{matrix} \quad t \quad
\\ l \,\, \boxe \,\, r  \\ \quad b\quad
\end{matrix}$ and $B = \begin{matrix} \quad u \quad \\ s \,\, \boxe \,\, m \\ \quad c \quad
\end{matrix}$ in $\B$.
\begin{flalign} \label{ax2-1} & \text{If } \quad A \vert B, \quad  \text{ then }  \quad AB =
\begin{matrix} \quad tu \quad \\ l \,\, \boxe \,\, m \\ \quad bc \quad
\end{matrix}, &
\\ \label{ax2-2} & \text{If }   \qquad \displaystyle\frac{A}{B},  \quad \text{ then }
\quad \begin{matrix}A\vspace{-4pt}\\B\end{matrix} =
\begin{matrix} \quad t \quad \\ ls \,\, \boxe \,\, rm \\ \quad c \quad
\end{matrix}.&
\end{flalign}

The notation $\begin{tabular}{p{0,4cm}|p{0,4cm}} $A$ & $B$ \\
\hline $C$ & $D$ \end{tabular}$ means that all possible horizontal
and vertical products are allowed; this implies that
$\displaystyle \frac{AB}{CD}$, $\begin{matrix} A \vspace{-4pt}\\ C\end{matrix} \Big\vert \begin{matrix} B \vspace{-4pt}\\ D\end{matrix}$.

\emph{Interchange law}. If
$\begin{tabular}{p{0,4cm}|p{0,4cm}} $A$ & $B$ \\ \hline $C$ & $D$
\end{tabular}$, then
\begin{equation}\label{permut}
\begin{matrix} A B \vspace{-4pt}\\ C D  \end{matrix} :=
\begin{matrix} \{A B\} \vspace{-2pt}\\ \{C D\} \end{matrix} =
\left\{\begin{matrix} A  \vspace{-4pt}\\ C   \end{matrix}\right\}
\left\{\begin{matrix} B \vspace{-4pt}\\ D  \end{matrix}\right\}.
\end{equation}

The identity functions $\idd: \Hc \to \B$ (vertical
identity), $\idd: \Vc \to \B$ (horizontal identity) satisfy $$
\idd(x) =
\begin{matrix} \quad x \quad \\ \begin{tabular}{||p{0,1cm}||} \hline \\
\hline \end{tabular} \\
\quad x \quad
\end{matrix}, \quad x \in \Hc; \qquad
\idd(g) = \begin{matrix}  g \,\, \begin{tabular}{|p{0,1cm}|} \hhline{|=|} \\
\hhline{|=|} \end{tabular} \,\, g
\end{matrix}, \quad g\in \Vc.
$$

We shall use the notation $\Theta_P: = \idd \, {\id_\Hc P} = \idd
\, {\id_\Vc P}$, $\forall P\in \Pc$.

\medbreak Suppose that
$A = \begin{matrix} \quad t \quad
\\ l \,\, \boxe \,\, r  \\ \quad b\quad
\end{matrix}$. The horizontal and vertical inverses of $A$ will be denoted by $A^h = \begin{matrix} \quad t^{-1} \quad
\\ r \,\, \boxee \,\, l  \\ \quad b^{-1}\quad \end{matrix}$, and
$A^v = \begin{matrix} \quad b \quad
\\ l^{-1} \,\, \boxe \,\, r^{-1}  \\ \quad t \quad \end{matrix}$, respectively.
The element $(A^h)^v = (A^v)^h$ will be denoted $A^{-1}$; thus
$A^{-1} = \begin{matrix} \quad b^{-1} \quad
\\ r^{-1}  \,\, \boxee \,\, l^{-1}   \\ \quad t^{-1}\quad \end{matrix}$.

If $\begin{tabular}{p{0,4cm}|p{0,4cm}} $A$ & $B$ \\ \hline $C$ &
$D$
\end{tabular}$, then
$$\left\{\begin{matrix} A B \vspace{-4pt}\\ C D  \end{matrix}\right\}^h
 = \begin{matrix} B^h \, A^h \vspace{-4pt}\\ D^h \,C^h  \end{matrix}, \quad \left\{\begin{matrix} A B \vspace{-4pt}\\ C D  \end{matrix}\right\}^v
 = \begin{matrix} C^v \, D^v \vspace{-4pt}\\ A^v \, B^v  \end{matrix}   \quad\text{and  }\left\{\begin{matrix} A B \vspace{-4pt}\\ C D  \end{matrix}\right\}^{-1}
 = \begin{matrix} D^{-1}\, C^{-1} \vspace{-4pt}\\ B^{-1}\, A^{-1}
\end{matrix}.$$

\bigbreak
\subsection{Some properties}

\

Let $\T$ be a double groupoid.
In this section we collect some technical results from
\cite[1.4]{AN} needed later.

\begin{lema}\label{l-unit}\cite[Lemma 1.9]{AN}.
Let $A, B, C \in \B$. The following statements are equivalent:

(i) $ABC \in \Hc$;

(ii) there exist $U, V \in \B$ such that
$\begin{matrix}U \vspace{-4pt}\\ V \end{matrix} = B$, $AU \in \Hc$, $VC \in \Hc$;

(iii) there exist $W, Z \in \B$ such that
$\begin{matrix}W \vspace{-4pt}\\ Z \end{matrix} = B$, $AZ \in \Hc$, $WC \in \Hc$.

Moreover, in (ii) and (iii) the elements $U, V, W, Z$ are
uniquely determined by $A, B, C$, and we have

$$\begin{tabular}{p{1,2cm}|p{0,4cm}|p{1,2cm}} $A$ & $U$ & $\iddv \,{t(C)}$ \\ \hline $\iddv \,{b(A)}$ & $V$ &
$C$ \end{tabular}, \qquad \text{respectively} \qquad
\begin{tabular}{p{1,2cm}|p{0,4cm}|p{1,2cm}} $\iddv \,{t(A)}$ & $W$ & C
\\ \hline A & $Z$ & $\iddv \,{b(C)}$ \end{tabular}.$$
\qed
\end{lema}

Dually, we have:

\begin{lema}\label{l-counit}\cite[Lemma 1.10]{AN}.
Let $A, B, C \in \B$. The following statements are equivalent:

(i) $\begin{matrix} A\vspace{-4pt}\\ B
\vspace{-4pt}\\ C \end{matrix} \in \Vc$;

(ii) there exist $U, V \in \B$ such that
$UV = B$, $\begin{matrix} A\vspace{-4pt}\\ U
\end{matrix} \in \Vc$, $\begin{matrix} V\vspace{-4pt}\\ C \end{matrix} \in \Vc$;

(iii) there exist $W, Z \in \B$ such that
$WZ = B$, $\begin{matrix} A\vspace{-4pt}\\ Z \end{matrix} \in \Vc$,
$\begin{matrix} W\vspace{-4pt}\\ C \end{matrix} \in \Vc$.

The elements $U, V, W, Z$ in (ii) and (iii) are
uniquely determined by $A$, $B$, $C$, and we have
$$\begin{tabular}{p{1,2cm}|p{1,2cm}} $A$ &  $\iddh \,{r(A)}$
\\ \hline $U$ & $V$ \\ \hline  $\iddh \,{l(C)}$ & $C$  \end{tabular},
\quad \text{respectively} \quad
\begin{tabular}{p{1,25cm}|p{1,25cm}} $\iddh \,{l(A)}$ & $A$
\\ \hline $W$ & $Z$ \\ \hline  $C$ &  $\iddh \,{r(C)}$ \end{tabular}.$$
\qed \end{lema}

\begin{lema}\label{l-atp3}\cite[Lemma 1.11]{AN}.
(i)  Let $A, X, Y, Z \in \B$ such that
\begin{equation} \label{tricot1}
\begin{tabular}{p{0,4cm}|p{0,8cm}|p{0,4cm}}\quad & $X^{-1}$
& \quad \\  \hline $X$ & $Y$ & $Z$ \\ \hline \quad  & $Z^{-1}$ & \quad  \end{tabular}\quad .
\end{equation}

Then the following conditions are equivalent:

\begin{align}
\label{tricot2}
XYZ &= A.
\\ \label{tricot25} \left\{ \begin{matrix} X^{-1} \vspace{-4pt}\\ Y \vspace{-4pt}\\ Z^{-1} \end{matrix} \right\} &= A^{-1}.
\end{align}

(ii) The collection $X = A = Z$, $Y = A^h$ satisfies
\eqref{tricot1}, \eqref{tricot2} and \eqref{tricot25}.

\medbreak
Moreover, if equation \eqref{tricot1} holds, then we have also
\begin{equation} \label{tricot15}
\begin{tabular}{p{0,7cm}|p{0,8cm}|p{0,7cm}}$X^v$ & $X^{-1}$ & $A^v$ \\
\hline $X$ & $Y$ & $Z$ \\ \hline $A^v$  & $Z^{-1}$ & $Z^v$  \end{tabular}\quad .
\end{equation}
\qed
\end{lema}

\bigbreak
\subsection{Core groupoids}\label{core}

\

Let $\T$ be a double groupoid.
Brown and Mackenzie have introduced a \emph{core} groupoid, and a
related core diagram, relevant to the structure of $\T$. Its
relation to the structure of a double Lie groupoid can be found in
\cite{mk2, BM}. In fact, there are four different core groupoids,
all isomorphic via the bijections given by vertical, horizontal and
total inversions. In this subsection we recall two of these core groupoids.
These will play an important r\^ ole in the
description of the source and target subalgebras of the quantum
groupoids constructed in Section \ref{wha}. Let

\begin{align*}\Dc
& : = \{ D \in \B: \; l(D), b(D) \in \Pc \},\\ \Ec & : = \{ E \in
\B: \; r(E), t(E) \in \Pc \}.\end{align*}

Thus elements of $\Dc$, resp. of $\Ec$, are of the form
\begin{tabular}{||p{0,1cm}|} \hline \\ \hhline{||=|}
\end{tabular}, resp. \begin{tabular}{|p{0,1cm}||} \hhline{|=||} \\
\hline
\end{tabular}.
Note that $\Theta_P \in
\Dc, \Ec$, for all $P \in \Pc$; so $\Dc, \Ec \neq \emptyset$.

\begin{prop} There are groupoid structures $s, e: \Dc \rightrightarrows
\Pc$, $s, e: \Ec \rightrightarrows \Pc$, with source and target
maps $s(D) = rt(D)$, $e(D) = lt(D)$, $D \in \Dc$, resp. $s(E) = bl(E)$,
$e(E) = br(E)$, $E \in \Ec$,  identity maps $\id: \Pc \to \Dc$, resp. $\id: \Pc \to \Ec$,  $P
\mapsto \Theta_P$, and compositions $\Dc \Times{e}{s} \Dc \to
\Dc$, $\Ec \Times{e}{s} \Ec \to \Ec$,  given by \begin{equation}D
\diamond L : = \left\{\begin{matrix} \iddv t(L)& D \vspace{-4pt}\\
L & \iddv r(L)\end{matrix} \right\}, \qquad M \circ E : =
\left\{\begin{matrix}\iddv l(E)& E \vspace{-4pt}\\ M &\iddv
b(E)\end{matrix} \right\}.
\end{equation} $D, L \in \Dc$, $M, E \in \Ec$.
The inverses of $D \in \Dc$ and $E \in \Ec$ are
\begin{align}\label{inversad}
D^{[-1]}: &= (\iddv t(D)^{-1}D)^v,
\\ \label{inversae}
E^{(-1)}: &= (E \iddv b(E)^{-1})^v
= \left\{\begin{matrix}\iddv l(E)^{-1} \vspace{-4pt}\\ E^h\end{matrix} \right\}.
\end{align}
 The map $D \mapsto D^{-1}$ gives an
isomorphism of groupoids $\Dc \overset{\simeq}\to \Ec$.\qed
\end{prop}

\begin{obs} $\Dc$ and $\Ec$ are \emph{not} subgroupoids of
$\B$; in particular, the inverses $D^{-1}$ and $D^{[-1]}$, $D \in
\Dc$, etc., should not be confused. Note that
\begin{equation}\label{daga}
D^{\daga}: = (D^{[-1]})^{-1} =  D^h \iddv t(D), \qquad D\in \Dc,
\end{equation}
defines an anti-isomorphism of groupoids, $(\quad)^{\daga}: \Dc \overset{\simeq}\to \Ec$, whose inverse is denoted by the same symbol.
\end{obs}

\begin{obs}\label{daga-appl}
If $DE = \iddh\, u$, for some $D\in \Dc$, $E\in \Ec$, $u\in \Hc$, then $u = t(D)$, $E = D^{\daga}$.
\end{obs}

\begin{obs} Recall that the restricted product of $\Hc^{\op}$ and $\Vc$
is the groupoid $$\Hc^{\op} \prode \Vc = \{(x, g) \in  \Hc \times
\Vc: l(x) = b(g),\ r(x) = t(g)\},$$ with componentwise
multiplication. Then there is a morphism of groupoids
$$\partial: \Dc \to
\Hc^{\op} \prode \Vc, \quad\partial(D)= (t(D), r(D)), \quad D\in
\Dc.$$

This morphism can be thought of as the \emph{core diagram}
of $\T$ introduced in \cite[Definition 2.1]{BM}; as shown in {\it
loc. cit.}, for all locally trivial double groupoids $\T$ (see Definition \ref{loctriv} below),
the core diagram of $\T$ determines $\T$.

The kernel of $\partial$ is the group bundle with $B\in \ker \partial(P)$
whenever its vertical and horizontal sides are the identities of
$P$, $P\in \Pc$.
\end{obs}

\medbreak There are several canonical maps from certain subsets of
$\B$ into $\Dc$ and $\Ec$. These arise naturally in the formulas
for the source and target maps for the weak Hopf algebras attached
to $\T$. They are described in what follows.

The formulas $$\phi(A) =
\left\{\begin{matrix} A^{-1} \vspace{-4pt}\\ \iddv
r(A)\end{matrix}\right\}, \qquad \alpha(A) =
\{ \iddv b(A)^{-1}A \},$$ define surjective maps
$$\phi: \{ A \in \B:\; t(A) \in \Pc \} \to \Dc, \quad
\alpha: \{ A \in \B:\; l(A) \in \Pc \} \to \Dc.$$
Note that if $D \in \Dc$, then $t(D^{-1}) \in \Pc$
and $D = \phi(D^{-1}) = \alpha(D)$.
On the other hand, the formulas $$\psi(A) = \left\{\begin{matrix}
\iddv l(A) \vspace{-4pt}\\ A^{-1}\end{matrix}\right\}, \qquad \beta(A)
= \{ A \iddv t(A)^{-1} \},$$ define  surjective maps
$$\psi: \{ A \in \B:\; b(A) \in \Pc \} \to \Ec, \quad
\beta: \{ A \in \B:\; r(A) \in \Pc \} \to \Ec.$$

\begin{prop}\label{acciones-core}
The groupoid $\Dc$ acts on the left on  the map $rt: \B \to
\Pc$, and on the right on  the map $rb: \B \to \Pc$, by the
formulas
$$D \rightharpoonup A = \left\{\begin{matrix} \iddv t(A)&
D \vspace{-4pt}\\ A & \iddv r(A)\end{matrix} \right\}, \qquad A
\leftharpoonup D = \left\{\begin{matrix}A & \iddv r(A) \vspace{-4pt}\\
\iddv b(A) & \iddv t(D)^{-1}D\end{matrix} \right\}.$$
Dually, $\Ec$ acts on the right on  the map $lt: \B \to \Pc$, and on the
left on the map $lb: \B \to \Pc$, by the formulas
$$A \leftharpoondown E = \left\{\begin{matrix} E\iddv b(E)^{-1}&
\iddv t(A) \vspace{-4pt}\\ \iddv l(A) & A\end{matrix} \right\},
\qquad E \rightharpoondown A = \left\{\begin{matrix} \iddv l(A)& A
\vspace{-4pt}\\ E & \iddv b(A)\end{matrix} \right\}.$$ \qed
\end{prop}
We remark for future use that there are bijections
\begin{align}\label{b-left}\{ A \in \B: \; l(A) \in \Pc \}
&\to \Hc \Times{r}{e} \Dc, \quad A \mapsto (b(A), \alpha(A)),
\\ \label{b-right}\{B \in \B: \; r(B) \in \Pc \}
&\to \Ec \Times{e}{l} \Hc, \quad B \mapsto (\beta(B), t(B)), \end{align}
the inverses given by horizontal composition.

\bigbreak
We next consider the action of the vertical composition groupoid
$\B \rightrightarrows \Hc$ on the map $\gamma: \Ec \to \Hc$, $\gamma(E) = b(E)^{-1}$, given by

\begin{equation}\label{curve}
A \curvearrowright E := \psi(E \rightharpoondown A)
= \left\{\begin{matrix}
\iddv \, l(A) \\  E\\ A^{-1} \end{matrix}\right\},
\end{equation}
$A\in \B$, $E\in \Ec$, $b(A) = b(E)^{-1}$; the second equality by
\eqref{inversae}.

We consider the following equivalence relation on the groupoid
$\Ec$: we say that $E$ and $M$ in $\Ec$ are \emph{vertically
connected}, denoted $E\sim_V M$ if, and only if, there exists
$g\in \Vc$ such that $t(g) = e(E)$, $b(g) = e(M)$; \emph{i. e.}
$g$ connects the ends of $E$ and $M$.

\begin{obs}\label{vertconn} Note that the following conditions are equivalent:

\begin{enumerate}
\item[(i)]  Any two elements of $\Ec$ are vertically connected;

\item[(ii)] $\Vc \rightrightarrows \Pc$ is
connected.\end{enumerate}

Indeed, the implication (ii) $\Rightarrow$ (i) is evident, while
(i) $\Rightarrow$ (ii) follows from vertical connectedness of the
boxes $\Theta_P \in \Ec$, $P \in \Pc$. \end{obs}

\begin{lema}\label{simvsiisimcurve}
Let $E$ and $M \in \Ec$. Then $E\sim_V M$ if and only
if there
exists $A\in \B$ such that $b(E)b(A) \in \Pc$ and
$M = A \curvearrowright E$.
\end{lema}

\pf If $M = A \curvearrowright E$, then $l(A) \in \Vc$ has source
$e(M)$ and target $e(E)$; hence
$E\sim_V M$. Conversely, assume that there exists $g\in \Vc$ such
that $t(g) = e(E)$, $b(g) = e(M)$. Then
$$
A = \left\{\begin{matrix} \vspace{-4pt} M^{-1}\\ \iddv g^{-1} \vspace{-4pt}\\
E^h\end{matrix}\right\}
$$
satisfies $M = \psi(E \rightharpoondown A) = = A \curvearrowright E$. \epf

\bigbreak
Finally, the maps $\phi$ and $\psi$ play also a r\^ole in the next result needed later.

\begin{lema}\label{conteo-st} Let $A \in \B$. Then we have

(i) There exist $X, Y \in \B$ such that
\begin{tiny}\begin{tabular}{p{0,25cm}|p{0,25cm}} \quad & $A$ \\
\hline $X$ & $Y$ \end{tabular}\end{tiny},
$XY \in \Hc$, $\begin{matrix}A \vspace{-4pt}\\ Y \end{matrix} \in \Vc
$, if and only if $t(A) \in \Pc$.
 In this case, we have  $Y = \phi(A)^h$ and $X = \iddv x \, \phi(A)$,
 for a unique $x \in \Hc$ such that $r(x) = rb(A)$.

(ii) There exist $X, Y \in \B$ such that
\begin{tiny}\begin{tabular}{p{0,3cm}|p{0,3cm}} $X$ & $Y$ \\
\hline $A$ & \quad \end{tabular}\end{tiny}, $XY \in \Hc$,
$\begin{matrix}X \vspace{-4pt}\\ A \end{matrix} \in \Vc$, if and
only if $b(A) \in \Pc$. In this case, we have $X = \psi(A)^h$ and
$Y = \psi(A)\, \iddv y$, for a unique $y \in \Hc$ such that $l(y)
= bl(A)$. \end{lema}

\pf (i) It is clear that if there is such a pair $X, Y$ then
$t(A)$ is an identity. Suppose that this is the case, and  let $X,
Y$ as in (i). Since $XY = \iddv t(XY)$,
then $r(Y) = \id tr(Y)$, and thus $\begin{matrix}A \vspace{-4pt}\\
Y \end{matrix} = \iddv r(A)$; on the other hand, since
$\begin{matrix}A \vspace{-4pt}\\ Y
\end{matrix}$ is an identity, then $b(Y) = \id lb(Y)$, and
therefore $X = \iddv b(X)Y^h = \iddv b(X)
\left\{\begin{matrix}A^{-1} \vspace{-4pt}\\ \iddv r(A) \end{matrix} \right\}$.
Note in addition that $rb(X) = br(A)$. Part (i) will be
established if we prove that for any $x \in \Hc$ with $r(x) =
br(A)$ there is an $X$ as in (i) with $b(X) = x$. This is done by letting
$X = \iddv x \left\{\begin{matrix}A^{-1} \vspace{-4pt}\\ \iddv r(A) \end{matrix} \right\}$.
We omit the proof of part (ii), which is similar. \epf

\bigbreak
\subsection{Corner functions}

\

We begin by introducing four 'corner' maps on the set of boxes. These will prove useful
later in order to define the coproduct of a related quantum groupoid.
Let
\begin{align*}& \ulcorner: \Vc {}_t\times_l\Hc \to \mathbb N \cup \{ 0 \},
\quad \llcorner: \Vc {}_b\times_l\Hc \to \mathbb N \cup \{ 0 \},
\\ & \urcorner: \Vc \Times{t}{r}\Hc \to \mathbb N \cup \{ 0 \},
\quad \lrcorner: \Vc {}_b\times_r\Hc \to \mathbb N \cup \{ 0
\},\end{align*} be given by the formulas:
\begin{align*}\ulcorner (g, x) & = \# \Big\{ U \in \B: \, U =  \begin{matrix} \quad \,\,\, x \quad \\  g \,\, \boxe \,\, \\ \quad \quad
\end{matrix} \Big\}; \quad
\urcorner (g, x) = \# \Big\{ U \in \B: \, U =
\begin{matrix} \quad x \,\,\, \quad \\  \,\, \boxe \,\, g \\ \quad \quad \end{matrix} \Big\};
\\
\llcorner (g, x) & = \# \Big\{ U \in \B: \, U =
\begin{matrix} \quad \quad \\  g \,\, \boxe \,\, \\ \quad \,\,\, x \quad \end{matrix} \Big\};
\quad \lrcorner (g, x) = \# \Big\{ U \in \B: \, U =
\begin{matrix} \quad \quad \\  \,\, \boxe \,\, g \\ \quad x \,\,\, \quad
\end{matrix} \Big\}.\end{align*}
Thus $\ulcorner (g, x)$ equals the number of boxes $B\in \B$ with
$l(B) = g$, $t(B) = x$, and so on.

\medbreak These determine four maps $\ulcorner, \llcorner,
\urcorner, \lrcorner: \B \to \mathbb N$ by the following rules:
\begin{align*}& \ulcorner (X) = \ulcorner (l(X), t(X)); \quad \urcorner (X) = \urcorner (r(X), t(X)); \\ & \llcorner (X) = \llcorner (l(X), b(X)); \quad \lrcorner (X) = \lrcorner (r(X), b(X)).\end{align*}

\begin{lema}\label{vh-1} (i) Let $x \in \Hc$, $g \in \Vc$. Then $\ulcorner(g, x) = \llcorner(g^{-1}, x) = \urcorner(g, x^{-1}) = \lrcorner(g^{-1}, x^{-1})$.

(ii) Let $X \in \B$. Then $\ulcorner(X) = \llcorner(X^v) = \urcorner(X^h) = \lrcorner(X^{-1})$.\end{lema}

\pf Part (i) is an easy consequence of the definitions, using the bijections given by vertical, horizontal and total inversions in $\T$.

Part (ii) follows from (i).
\epf

\begin{definition} Let $A, B, X, Y \in \B$ such that
$XY = \begin{matrix} A \vspace{-4pt}\\ B \end{matrix}$. By a \emph{double factorization} of this common product we shall mean a quadruple $(U, V, R, S)$ of elements in $\B$ satisfing
\begin{equation}\label{cuadrado}\begin{tabular}{p{0,4cm}|p{0,4cm}} $U$ &
$V$ \\ \hline $R$ & $S$ \end{tabular}, \quad UV = A, \quad RS = B,
\quad \begin{matrix} U \vspace{-4pt}\\ R \end{matrix}= X, \quad \begin{matrix} V \vspace{-4pt}\\
S \end{matrix} = Y.\end{equation}
The set of all double factorizations will be denoted $[X, Y, A, B]$. \end{definition}

\begin{prop}\label{rel-corners} Let $A, B, X, Y \in \B$ such that $XY = \begin{matrix} A \vspace{-4pt}\\ B \end{matrix}$.  Then we have
\begin{align*}\#[X, Y, A, B] & = \ulcorner(l(A), t(X)) =  \llcorner(l(B), b(X)) \\ & =
\urcorner(r(A), t(Y)) =  \lrcorner(r(B), b(Y)).\end{align*} \end{prop}

\pf The map $[X, Y, A, B] \to \Big\{ U \in \B: \, U =  \begin{tiny}\begin{matrix} \quad \qquad t(X) \quad \\  l(A) \,\, \boxe  \\ \quad \quad \end{matrix}\end{tiny} \Big\}$, given by $(U, V, R, S) \mapsto U$ is a well defined bijection, whose inverse is given by  $U \mapsto \Big(U, U^hA, \begin{matrix} \; U^v \vspace{-4pt}\\ X \end{matrix}, \; \begin{matrix} U^{-1}A^v \vspace{-4pt}\\ Y \end{matrix}\Big)$.  This shows the first equality. The others are similarly established. \epf

As a consequence we get the following symmetry properties of the
corner maps:

\begin{cor}\label{simetria} Let $L, M, N \in \B$. Suppose that
\begin{tiny}\begin{tabular}{p{0,25cm}|p{0,3cm}} $L$ &
$M$ \\ \hline $N$ & \quad \end{tabular}\end{tiny}. Then we have
\begin{align*} (i)& \quad \ulcorner(L) = \urcorner(M), \qquad
(ii)\quad \llcorner(L) = \lrcorner(M), \\  (iii)& \quad
\ulcorner(L) = \llcorner(N), \qquad (iv)\quad \urcorner(L) =
\lrcorner(N).\end{align*} \end{cor}

In particular, $\ulcorner(L) = \urcorner(L)$, $\llcorner(L) = \lrcorner(L)$, for every horizontal identity $L \in \Vc$, and $\ulcorner(N) = \llcorner(N)$, $\urcorner(N) = \lrcorner(N)$, for every vertical identity $N \in \Hc$.

\pf The quadruple $X = L$, $Y = M$, $A = LM$, $B = \id b(LM)$
satisfies the assumptions of Proposition \ref{rel-corners}, and
then the proposition implies that $$\ulcorner(l(L), t(L)) =
\ulcorner(l(A), t(X)) =  \urcorner(r(A), t(Y)) = \urcorner(r(M),
t(M));$$ this proves part (i). Part (ii) follows similarly,
considering instead the set $[L, M, \iddv t(LM), LM]$. As to parts
(iii) and  (iv),  the same arguments apply with $A = L$, $B = N$,
$X = \begin{matrix} L \vspace{-4pt}\\ N \end{matrix}$, $Y = \id
r\left(\begin{matrix} L \vspace{-4pt}\\ N \end{matrix}\right)$.
\epf

\begin{lema}\label{corner-p} Let $P \in \Pc$.
Then we have
\begin{equation*}\lrcorner(\id_{\Vc} P, \id_{\Hc} P)
= \urcorner(\id_{\Vc} P, \id_{\Hc} P) = \llcorner(\id_{\Vc} P, \id_{\Hc} P)
= \ulcorner(\id_{\Vc} P, \id_{\Hc} P). \end{equation*} \end{lema}

\emph{The common value in Lemma \ref{corner-p} will be denoted
$\theta(P)$.} This agrees with the value of any of the corner
functions on the box $\Theta_P$.

\pf The proof follows from Lemma \ref{vh-1}.\epf

A surprising consequence of Proposition \ref{rel-corners} is that the corner functions on a box actually depend only on one vertex, the vertex 'opposite' to the corner, of that box.

\begin{prop}\label{puntos} Let $L \in \B$. Then the following hold.
\begin{align*}(i)&  \; \ulcorner(L) = \theta(br(L)), \qquad
(ii)\; \llcorner(L) = \theta(rt(L)),   \\ (iii)& \;
\urcorner(L) = \theta(bl(L)), \qquad (iv)\; \lrcorner(L) = \theta(tl(L)).\end{align*} \end{prop}

\pf We show part (i); then parts (ii)--(iv) follow from (i) and
Lemma \ref{vh-1}. To do this we argue as in the proof of Corollary
\ref{simetria}. Let $X = L$, $Y = \iddv r(L)$, $A = L$, $B = \iddv
b(L)$. By Proposition \ref{rel-corners}, $$\ulcorner(L) =
\ulcorner(l(A), t(X)) = \lrcorner(r(B), b(Y)) =
\lrcorner(id_{\Vc}br(L), \id_{\Hc} br(L)) = \theta(br(L)),$$ as
claimed.  \epf

\begin{lema}\label{d-cnt} Let $P, Q \in \Pc$. Suppose that
$P$ and $Q$ are connected by an arrow of the core groupoid $\Dc$.
Then $\theta(P) = \theta(Q)$. \end{lema}

Note that $P$ and $Q$ are connected by $\Dc$ if and only if they are connected by $\Ec$.

\pf Let $[P]$ denote the connected component of $\Pc$  with respect to $\Dc$ containing $P$. We have
\begin{align*}\theta(P) & = \sum_{R \in \Pc} \# \{ B \in \B: \; l(B) = \id_P, \, b(B) = \id_P, \, rt(B) = R \} \\ & = \sum_{R \in \Pc} \# \Dc(R, P) = \Dc(P)\; \#[P].\end{align*} Thus, $\theta(P) = \theta(Q)$, whenever they are in the same connected component with respect to $\Dc$.\epf

The following proposition states  the main translation invariance
properties of the above defined maps.

\begin{prop}\label{trans-inv} Let $X, Y, Z \in \B$ such that \begin{tiny}\begin{tabular}{p{0,25cm}|p{0,25cm}} $X$ &
$Y$ \\ \hline $Z$ & \quad \end{tabular}\end{tiny}. Then we have
\begin{align*}
(i)& \quad \urcorner(XY) = \urcorner(X), \qquad (ii)\quad
\urcorner\left(\begin{matrix} X \vspace{-4pt}\\ Z \end{matrix}\right) =
\urcorner(Z),\\ (iii)& \quad \llcorner(XY) = \llcorner(Y),\qquad
(iv)\quad \llcorner\left(\begin{matrix} X \vspace{-4pt}\\ Z \end{matrix}\right)
= \llcorner(X).\end{align*}
\end{prop}

Similar properties hold for the functions $\lrcorner$ and $\ulcorner$.

\pf  We prove (i): $$\urcorner(XY) =\theta(bl(XY)) =\theta(bl(X))=  \urcorner(X),$$
by
Proposition \ref{puntos}. The proof of
parts (ii), (iii) and (iv) is similar.  \epf

\bigbreak
\subsection{Counting formulas in double groupoids}

\

We obtain in this subsection some counting formulas, involving the corner functions, that will be of use in the next section.

\begin{lema}\label{conteo-ant} Let $A \in \B$. There is a bijection between $$\big\{ (X, Y, Z) \in \B^3: \, (X, Y, Z) \,
\text{satisfies} \, \eqref{tricot1}, \eqref{tricot2} \, \text{in }
\, \ref{l-atp3} \big\}$$ and $\Big\{X \in \B: X =  \begin{tiny}
\begin{matrix} \quad \,\, \quad \\  l(A) \boxe \,\, \quad \\ \quad b(A)\end{matrix}\end{tiny} \Big\} \times \Big\{Z \in \B: \, Z =  \begin{tiny} \begin{matrix} \, t(A) \quad \quad \\  \,\,
\boxe \,\, r(A) \\ \quad \quad \end{matrix}\end{tiny} \Big\}$. In particular,
\begin{equation*}\#\big\{ (X, Y, Z) \in \B^3: \, (X, Y, Z) \,
\text{satisfies} \, \eqref{tricot1}, \eqref{tricot2} \, \text{in }
\, \ref{l-atp3} \big\} = \llcorner(A)\urcorner(A).
\end{equation*}
\end{lema}

\pf For any such triple we have $X = $\, \begin{tiny}
$\begin{matrix} \quad \,\, \quad \\  l(A) \boxe \,\, \quad \\ \quad b(A)\end{matrix}$\end{tiny}
and $Z = $  \begin{tiny}$\begin{matrix} \, t(A) \quad \quad \\  \,\,
\boxe \,\, r(A) \\ \quad \quad \end{matrix}$\end{tiny}.
Moreover, $Y$ is determined by $Y = X^hAZ^h$. Conversely, for
every pair $(X, Z)$ as above, the triple $(X, X^hAZ^h, Z)$ satisfies \eqref{tricot1}, \eqref{tricot2}. This proves the lemma. \epf

\begin{lema}\label{u-v} Let $A, X, Y \in \B$. Assume that $tr(A) = lb(Y)$ and $bl(A) = rt(X)$.
Then the following hold:

\begin{multline}\label{i} \# \big\{ (U, V) \in \B^2: \,
 UV = A,
\, \begin{tiny}\begin{matrix} U \\ \quad
V^{-1}\end{matrix}\end{tiny} = Y  \big\}
\\ = \begin{cases}\urcorner(r(A),
b(Y)^{-1})  = \ulcorner(l(A), t(Y)) \quad \text{if } b(A) \in \Pc,
\\0 \qquad\qquad\qquad\qquad\qquad\qquad\qquad \text{otherwise.}\end{cases}
\end{multline}
\begin{multline}\label{ii} \# \big\{ (U, V) \in \B^2: \, UV = A, \,
\begin{tiny}\begin{matrix} \quad U^{-1} \\ V \end{matrix}\end{tiny}
= X  \big\}  \\ = \begin{cases}\llcorner(l(A), t(X)^{-1}) = \lrcorner(r(A),
b(X))
\quad \text{if } t(A) \in \Pc,
\\0 \qquad\qquad\qquad\qquad\qquad\qquad\qquad \text{otherwise.}\end{cases}
\end{multline} \end{lema}

\pf The set of pairs $(U, V)$ as in \eqref{i} coincides with the
set of pairs $(U, V)$ satisfing
\begin{equation*}\begin{tabular}{p{0,7cm}|p{0,7cm}} $U$ &
$V$ \\ \hline $V^{-1}$ & $V^v$ \end{tabular}, \quad UV = A, \quad
\begin{matrix} U \vspace{-4pt}\\ V^{-1}\end{matrix} = Y, \end{equation*} which
is in bijective correspondence with the set of quadruples $(U, V,
R, S)$ such that
\begin{equation*}\begin{tabular}{p{0,4cm}|p{0,4cm}} $U$ &
$V$ \\ \hline $R$ & $S$ \end{tabular}, \quad UV = A, \quad
\begin{matrix} U \vspace{-4pt}\\ R\end{matrix} = Y, \quad RS = \iddv r(A)^{-1},
\quad \begin{matrix} V \vspace{-4pt}\\ S\end{matrix} = \iddv b(Y)^{-1},
\end{equation*} by sending $(U, V)$ to $(U, V, V^{-1}, V^v)$, and
$(U, V, R, S)$ to $(U, V)$. Hence Equation \eqref{i} follows from
Proposition \ref{rel-corners}. Equation \eqref{ii} is similarly
shown. \epf

\bigbreak
\subsection{Vacant double groupoids}

\

A double groupoid $\T$ is called \emph{vacant} if for any $g\in \Vc$, $x\in \Hc$
such that $r(x) = t(g)$, there is exactly one  $X \in \B$ such that
$X = \begin{matrix} \quad x \quad \\  \,\, \boxe \,\, g \vspace{-6pt}\\ \quad \quad \end{matrix}$.
Vacant double groupoids have been introduced in \cite[Definition 2.11]{mk1}. We gave in \cite{AN} several characterizations of vacant double groupoids that we
 had found in the course of our research; the following is an example of a completely
 symmetric one. Its proof is a direct consequence of Proposition \ref{rel-corners}.

\begin{prop}\label{vert-hor}\textbf{\cite{AN}}. Let $\T$ be a double groupoid.
The following are equivalent.
\begin{enumerate}
\item $\T$ is vacant. \item For all $A, B, X, Y \in \B$ such that
$XY = \begin{matrix} A \vspace{-4pt}\\ B \end{matrix}$, there exist unique $U$,
$V$, $R$, $S \in \B$ satisfing Equation \ref{cuadrado}.
\end{enumerate} \qed\end{prop}

In terms of corner functions, vacant double groupoids are
characterized by the property that some (hence all) corner
function $\Vc {}_{t/b}\times_{l/r}\Hc \to \mathbb N \cup \{ 0 \}$
takes constantly the value $1$.

\medbreak If $\T$ is vacant, the core groupoids $\Dc$ and $\Ec$
are isomorphic and coincide with the discrete groupoid on the base
$\Pc$: that is, the only arrows in the core groupoids are the
identity arrows. In fact, we have the following characterization.

\begin{prop}\label{vacant-corner} Let $\T$ be a double groupoid.
The following are equivalent.
\begin{enumerate}
\item $\T$ is vacant.
\item Some (hence all) corner function takes positive values\footnote{This is the filling condition \eqref{neq0} below.} and the core groupoids are discrete on the base
$\Pc$.
\end{enumerate} \end{prop}

\pf We have already discussed $1\implies 2$.

$2\implies 1$. Since $\Ec \simeq \Pc$, $\theta(P) = 1$ for all $P\in \Pc$. If $(g,x)\in \Vc {}_{t}\times_{r}\Hc$ then there exists at least one box $A\in \B$ with $t(A) = x$, $r(A) = g$. Let $P = bl(A)$. Then $\urcorner(g,x) = \urcorner(A) = \theta (P) = 1$, the second equality by Proposition \ref{puntos} (iii).
\epf

\section{Construction of quantum groupoids from double groupoids}\label{wha}

\bigbreak
\subsection{Quantum groupoids}

\

Recall \cite{bnsz, bsz} that a \emph{weak bialgebra} structure on a
vector space $H$ over a field $\ku$ consists of an associative
algebra structure $(H, m, 1)$, a coassociative coalgebra structure
$(H, \Delta, \epsilon)$, such that the following are satisfied:
\begin{align}\label{d-mult}\Delta(ab) &= \Delta(a) \Delta(b), \qquad \forall a, b \in H.
\\ \label{ax-unit} \Delta^{(2)} (1 ) &= \left( \Delta(1) \otimes 1 \right)
\left( 1 \otimes \Delta(1) \right) = \left( 1 \otimes \Delta(1)
\right) \left( \Delta(1) \otimes 1 \right).
\\ \label{ax-counit} \epsilon(abc) &= \epsilon(ab_1)\epsilon(b_2c)
= \epsilon(ab_2)\epsilon(b_1c), \qquad \forall a, b, c \in
H. \end{align}
The maps $\epsilon_s$, $\epsilon_t$ given by
\begin{align*}  \epsilon_s(h) &
= (\id \otimes \epsilon) \left( (1 \otimes h) \Delta(1) \right),
\\\epsilon_t(h) &
= (\epsilon \otimes \id) \left(  \Delta(1) (h \otimes 1)\right),
\end{align*}
$h \in H$, are
respectively called the source and target maps; their images are
respectively called the source and target subalgebras.

A weak bialgebra $H$ is called a \emph{weak Hopf algebra} or a
\emph{quantum groupoid} if there exists a linear map $\Ss: H \to
H$ satisfying
\begin{align}\label{atp-1} m(\id \otimes \Ss) \Delta (h) &
=  \epsilon_t(h), \\
\label{atp-2} m(\Ss \otimes \id) \Delta (h) &
= \epsilon_s(h),\\
\label{atp-3} m^{(2)}(\Ss \otimes \id \otimes \Ss) \Delta^{(2)} &
= \Ss,
\end{align}
for all $h \in H$.  See
\cite{nik-v} for a survey on quantum groupoids. It is known that a
weak Hopf algebra is a true Hopf algebra if and only if $\Delta(1)
= 1 \otimes 1$.

\bigbreak
\subsection{Weak Hopf algebras arising from double groupoids}

\

Let $\T$ be a  \emph{finite} double groupoid,  that is,  $\B$ is a finite set (and \emph{a fortiori} also $\Vc$, $\Hc$ and $\Pc$ are finite).

Let $\ku$ be a field of characteristic zero and let $\ku\T$ denote
the $\ku$-vector space with basis $\B$. We define a multiplication
and comultiplication on $\ku\T$ by the formulas
\begin{equation}\label{prod}A.B = \begin{cases} \begin{matrix}A \vspace{-4pt}\\B\end{matrix}, \quad
\text{if } \displaystyle\frac{A}{B}, \\
0 , \quad \text{otherwise}, \end{cases} \end{equation}
\begin{equation}\label{del} \Delta(A) = \sum_{XY = A} \dfrac{1}{\urcorner(Y)} \; X\otimes Y
= \sum_{XY = A} \dfrac{1}{\ulcorner(X)} \; X\otimes Y,
\end{equation} for all  $A, B \in \B$. The second identity because of Corollary \ref{simetria}.

\medbreak Therefore $\ku \T$ is an associative algebra with unit
$\uno : = \sum_{x \in \Hc} \iddv \,{x}$. This algebra structure
coincides with the groupoid algebra structure on $\ku\T$
corresponding to the vertical composition groupoid $\B \rightrightarrows \Hc$. See
\cite{AN}.

The coalgebra structure
on $\ku \T$ is a modification of the dual  groupoid coalgebra of
the horizontal composition  groupoid $\B \rightrightarrows \Vc$, studied in \cite{AN}.

\begin{lema} The comultiplication \eqref{del} makes $\ku \T$ into a
coassociative coalgebra, with counit $\epsilon: \ku\T \to \ku$
given by $$\epsilon(A) = \begin{cases} \ulcorner(A) =\urcorner(A), \quad
\text{if } A \in \Vc,  \\ 0 , \quad \text{otherwise}.
\end{cases}$$ \end{lema}

\pf Let $A \in \B$. We have
\begin{align*}(\Delta \otimes \id)\Delta(A) & = \sum_{XYZ = A}\dfrac{1}{\urcorner(Y) \urcorner (Z)} \; X \otimes Y \otimes
Z, \\ (\id \otimes \Delta)\Delta(A) & = \sum_{XYZ =
A}\dfrac{1}{\urcorner(YZ) \urcorner (Z)} \; X \otimes Y \otimes Z.
\end{align*} Thus coassociativity of $\Delta$ follows from
Proposition \ref{trans-inv}. The counit axiom $(\id \otimes
\epsilon)\Delta = \id$ is straightforward to check. Also, using
the definitions of $\Delta$ and $\epsilon$, we have for every $A
\in \B$, \begin{equation*}(\epsilon \otimes \id)\Delta(A) =
\sum_{XY = A} \dfrac{1}{\ulcorner(X)} \, \epsilon(X)Y =
\dfrac{\urcorner(l(A))}{\ulcorner(l(A))}\, A = A;
\end{equation*}
the last identity in view of Corollary \ref{simetria}.
\epf

\medbreak Recall from Subsection \ref{core} the definition of the
core groupoids $\Dc$ and $\Ec$. Let $D \in \Dc$, $E \in \Ec$. We
introduce elements ${}_D\uno, \uno_E \in \ku \T$ by
\begin{equation}{}_D\uno := \sum_{z \in \Hc, \, r(z) = e(D)} \{\iddv z \, D\},
\quad \uno_E := \sum_{x \in \Hc, \, l(x) = e(E)} \{E \, \iddv x\}. \end{equation} Observe that the maps $D \mapsto {}_D\uno$, $E \mapsto \uno_E$, are both injective.

\medbreak We relate the groupoid structures of $\Dc$, $\Ec$ with the
multiplication in $\ku \T$.

\begin{lema}\label{gpd-st} Let $D, L \in \Dc$, $E, M \in \Ec$. We have
\begin{align*}{}_D\uno . {}_L\uno &=
\begin{cases}{}_{D \diamond L}\uno, \; \text{if} \; e(D) = s(L),
\\ 0, \quad \text{otherwise};\end{cases}
\\ \uno_E . \uno_M &= \begin{cases}\uno_{M \circ E}, \; \text{if} \;
e(M) = s(E),
\\ 0, \quad \text{otherwise}.\end{cases}\end{align*}
\end{lema}

\pf We compute

\begin{align*}{}_D\uno . {}_L\uno & =
\sum_{\begin{tiny}\begin{matrix}r(z) = e(D), \\ r(w) = e(L)\end{matrix}\end{tiny}}
\{\iddv z \, D\}.\{\iddv w \, L\}
= \sum_{\begin{tiny}\begin{matrix}r(z)  = e(D), \\ r(w)  = e(L),
\\ z  = wt(L)\end{matrix}\end{tiny}}
\left\{\begin{matrix}\iddv z \, D  \vspace{-4pt}\\ \iddv w \, L\end{matrix}\right\}
\\
& = \delta_{e(D), s(L)} \sum_{r(w)  = e(L)}
\left\{\begin{matrix}\iddv w \, \iddv t(L) \, D
\vspace{-4pt}\\ \iddv w \, L\end{matrix}\right\}
\\ & = \delta_{e(D), s(L)} \sum_{r(w)  = e(L)} \iddv w \, (D \diamond L)
= \delta_{e(D), s(L)} \, {}_{D \diamond L}\uno.
\end{align*} This proves the first claim. The second is similarly shown.  \epf

\medbreak In what follows we shall consider double groupoids $\T$
which satisfy the following \emph{filling condition}:
\begin{equation}\label{neq0}\urcorner(g, x) \neq 0, \qquad
\forall x \in \Hc, \, g \in \Vc, \, r(x) = t(g).
\end{equation}This filling condition on double groupoids has been considered by Mackenzie \cite{mk2}. It is easy to exhibit examples of double groupoids that do not satisfy the filling condition;
\emph{e.~g.} the union of the vertical and horizontal subgroupoids of a suitable double groupoid.
There are three other equivalent formulations of \eqref{neq0} in terms of the other corners, \emph{cf.} Lemma \ref{vh-1}.

\begin{theorem}\label{bicross} Suppose that $\T$ is a double
groupoid satisfying \eqref{neq0}. Then $\ku\T$ is a weak Hopf
algebra with multiplication and comultiplication given by
\eqref{prod} and \eqref{del}. The antipode is determined by the
formula \begin{equation}\label{antipode}\mathcal S(A) =
\dfrac{\ulcorner(A)}{\llcorner(A)}\; A^{-1},\end{equation} for all
$A \in \B$. The source and target maps are given, respectively, by
\begin{align} \label{source} \epsilon_s(A) & =
\begin{cases}{}_{\phi(A)}\uno, \quad \text{if} \quad t(A) \in \Pc, \\
0, \qquad \text{otherwise};  \end{cases}  \\
\label{target} \epsilon_t(A) & =
\begin{cases}\dfrac{\ulcorner(A)}{\llcorner(A)}\; \uno_{\psi(A)}, \quad \text{if}
\quad b(A)  \in \Pc, \\
0, \qquad \text{otherwise}.  \end{cases}  \end{align} The source
and target subalgebras are isomorphic to the groupoid algebras
$\ku \Dc$ and  $(\ku \Ec)^{\op}$, respectively. \end{theorem}

Here, the maps $\phi, \psi$ are those defined in Subsection \ref{core}.

\begin{rmk} (i). The weak Hopf algebra $\ku\T$ is semisimple, since the underlying algebra is a groupoid algebra.

(ii). The map $\lambda: \B \to \ku^{\times}$, $\lambda (A) =
\dfrac{\ulcorner(A)}{\llcorner(A)}$ is a character with respect to vertical composition, see Lemma \ref{char}.

\end{rmk}

\medbreak Theorem \ref{bicross} will be proved in Section \ref{extcocycles}.
We find appropriate to make
the following observation. The algebra and coalgebra structures
depend intrinsically on the 'double' nature of operations in $\T$.
It would seem that the construction of $\ku \T$ shows some
preference for the vertical composition over the horizontal one.
However, as the following remark shows, there is no such
preference: the symmetry of the construction is hidden behind the
choice of basis in $\ku \T$.

\begin{obs}\label{nopref}Let $H$ be the vector space with basis
$\{\underline{A}\}_{A \in \B}$, and multiplication and
comultiplication defined by
\begin{align}\underline{A}.\underline{B} &=
\begin{cases} \dfrac{1}{\urcorner(A)}\;\underline{C}
= \dfrac{1}{\lrcorner(B)}\;\underline{C},
\quad \text{if } \displaystyle\frac{A}{B}, \quad
\text{where} \; C = \begin{matrix}A \vspace{-8pt}\\B\end{matrix};\\
0 , \quad \text{otherwise}, \end{cases} \\
\Delta(\underline{A}) &= \sum_{XY = A}  \underline{X}\otimes \underline{Y},
\end{align} for all  $A, B \in \B$.
Then the linear isomorphism $H \to \ku \T$, $\underline{A} \mapsto
\dfrac{1}{\urcorner(A)}A$, preserves multiplication and
comultiplication, in view of Proposition \ref{trans-inv} (ii) and
(i), respectively. \end{obs}

\bigbreak
\subsection{The square of the antipode}

\

In this subsection we give the relations between corner functions
and the square of the antipode. The  proof of the following lemma
is straightforward.

\begin{lema}\label{sq-antp} Let $A \in \B$. Then $\mathcal S^2(A) =
\frac{\ulcorner(A)\lrcorner(A)}{\llcorner(A)\urcorner(A)}A$. \qed
\end{lema}

It follows from Proposition \ref{puntos}, that if $tl(A) = bl(A)$
and $br(A) = tr(A)$, then $\mathcal S^2(A) = A$.

\bigbreak
Recall from \cite[Definition 2.2.2]{nik-ss} that a weak Hopf algebra is called \emph{regular} if $\mathcal S^2 = \id$
on the source and target subalgebras.  Actually, if $\mathcal S^2 = \id$ on
the source (respectively, target) subalgebra, then also $\mathcal S^2 = \id$ on
the target (respectively, source) subalgebra.

\begin{prop}\label{regular} $\ku \T$ is a regular weak Hopf algebra.
\end{prop}

\pf Let $D \in \Dc$ and $E \in \Ec$. Then
$$
\Ss(\uno_E) = {}_{E^{-1}}\uno, \qquad \Ss({}_D\uno) =
\frac{\theta(s(D))}{\theta(e(D))}\, \uno_{D^{-1}} = \uno_{D^{-1}},
$$
the last equality by Lemma \ref{d-cnt}.
Hence $\mathcal S^2 = \id$ on the
source subalgebra. This proves the proposition. \epf

\bigbreak
It is natural to ask whether $\Ss^2 = \id$ on the whole weak
Hopf algebra $\ku\T$. It is possible to give a positive answer in some cases, \emph{e. g.} if $\T$ is a double group (that is, if $\Pc$ has just one element), if $\T$ is vacant, or if $\T$ is locally trivial-- see Example
\ref{loctriv} below.
But there are double groupoids $\T$ with $\Ss^2 \neq \id$.

\begin{exa}\label{no-siempre} We show that, for any positive rational number $r$, there exists a double groupoid $\T$ satisfying \eqref{neq0}
and a box $A\in \B$ such that $\Ss^2(A) = rA$.

\bigbreak
Let $m, n$ be natural numbers. Let $\Pc$ be the set with $m+n+3$ elements labeled
$$
P, \,Q,\, R,\, S_1, \dots, \, S_m, \, T_1, \dots,\, T_n.
$$

Let $\Hc$ and $\Vc$ be the groupoids with base $\Pc$ corresponding to the equivalence relations $\sim_{\Hc}$ and $\sim_{\Vc}$
given, respectively, by the following partitions:

\begin{align*} \sim_{\Hc}
: \quad &  \{P, \,Q,\, T_1, \dots,\, T_n\}
\, \bigcup \, \{R,\, S_1, \dots, \, S_m\},
\\ \sim_{\Vc}: \quad
&  \{P,  R\}
\, \bigcup \,\{Q,\, S_1, \dots, \, S_m, \, T_1, \dots,\, T_n\}.
\end{align*}

Let $\T = \begin{matrix} \B &\rightrightarrows &\Hc
\\\downdownarrows &&\downdownarrows \\ \Vc &\rightrightarrows &\Pc \end{matrix}$ be the double groupoid consisting of commuting squares
in the coarse groupoid $\Pc \times\Pc$, with horizontal arrows
in $\Hc$ and vertical arrows in $\Vc$. We claim that $\T$ satisfies the filling condition \eqref{neq0}. Moreover, any $u\in \Pc
\times \Pc$ can be expressed as a product $gx$ with $g\in \Vc$,
$x\in \Hc$. Indeed, if $u\in \Hc$ or $\Vc$, this is clear. The remaining possibilities are:

\begin{align*}
u &= (P, S_j) = (P,R)(R, S_j), \\
u &= (Q, R) = (Q, S_j)(S_j,R), \\
u &= (S_j, P) = (S_j,Q)(Q, P), \\
u &= (R, Q) = (R,P)(P, Q), \\
u &= (R, T_i) = (R,P)(P, T_i), \\
u &= (T_i, R) = (T_i, S_j)(S_j, R).
\end{align*}

Next, if $X\in \Pc$ then $\theta(X)$ equals the number of arrows in $\Hc \cap \Vc$ with source $X$. In the present example, this number equals $$\# \{ Y \in \Pc:\, Y\sim_{\Vc}X, \, Y\sim_{\Hc}X \}.$$
Thus
$$
\theta(P) =1 = \theta(R), \qquad \theta(Q) = n + 1 = \theta(T_i),
\qquad \theta(S_j) = m,
$$
$1\le j\le m$, $1\le i \le n$. Consider now the box
$$
A = \begin{matrix}\quad P  \quad Q \quad
\\  \,\, \boxe  \,\,  \\ \quad R  \quad S_1 \quad \end{matrix}  \in \B.
$$
Then
$$
\Ss^2(A) = \frac{\ulcorner(A)\lrcorner(A)}{\llcorner(A)\urcorner(A)}A
= \frac{\theta(S_1)\theta(P)}{\theta(Q)\theta(R)}A =
\frac{m}{n + 1}A.
$$
\end{exa}

\bigbreak
\begin{exa}\label{loctriv} Let us first recall the following definitions  \cite[Definition 2.3]{BM}. Let $\T$ be a
finite double groupoid.

\begin{enumerate}
\item[(a)] $\T$ is \emph{horizontally transitive} if every configuration of matching sides
$$
\begin{matrix}  l \,\, \boxu \,\, r \\ \quad b\quad
\end{matrix},
$$
$l,r\in \Vc$, $b\in \Hc$, can be completed to at least one box in $\B$. Equivalently, with
$\boxp$ instead of $\boxu$.

\medbreak
\item[(b)] $\T$ is \emph{vertically transitive} if every configuration of matching sides
$$
\boxa
$$
can be completed to at least one box in $\B$. Equivalently, with
$\boxc$ instead of $\boxa$.

\medbreak
\item[(c)] $\T$ is \emph{transitive} or \emph{locally trivial}
if it is both vertically  and horizontally transitive.
\end{enumerate}

\bigbreak
Note that any double groupoid that is either horizontally or vertically transitive satisfies the filling condition \eqref{neq0}.
Indeed, assume that $\T$ is vertically transitive and fix
$(x,g) \in \Hc\Times{r}{t} \Vc$. Then the configuration $$
\begin{matrix} x  \quad \\
\boxa \quad g \\ \id b(g)  \quad
\end{matrix}$$
can be completed to at least one box in $\B$.

\medbreak
\begin{lema}\label{loctrivlema}
(i). If  $\T$ is vertically transitive then the corresponding
corner function $\theta$ is constant along the connected components of
$\Hc \rightrightarrows \Pc$. Hence $\Ss^2 = \id$ on $\ku\T$.

\medbreak
(ii). If  $\T$ is horizontally transitive then the corresponding
$\theta$ is constant along the connected components of
$\Vc \rightrightarrows \Pc$. Hence $\Ss^2 = \id$ on $\ku\T$.
\end{lema}

\pf (i). Let $P, Q\in \Pc$ and $x\in \Hc$ such that $P = l(x)$, $Q = r(x)$. Complete the configuration of matched sides
$$
\begin{matrix} \qquad x \\  \id P \,\, \begin{tabular}{||p{0,1cm}} \hline \\ \hhline{||=|}
\end{tabular} \\ \qquad\qquad \id P\quad
\end{matrix}
$$
to a box $D \in \Dc$. Then $D$ connects $P$ and $Q$, so $\theta(P) = \theta(Q)$ by Lemma \ref{d-cnt}.
The second claim follows from Lemma \ref{sq-antp}.

The proof of (ii) is analogous, or else follows from (i) by passing to the transpose double groupoid. \epf

\end{exa}

\bigbreak
\subsection{The category $\Rep \ku \T$}

\

Let $\T$ be a finite double
groupoid satisfying the filling condition \eqref{neq0}. Let $\ku \T$
be the quantum groupoid associated to $\T$ as in Theorem
\ref{bicross}. Consider the category $\C : = \Rep \ku \T$ of
finite dimensional representations of $\ku\T$. This is a
semisimple category, with a finite number of simple objects.

Objects of $\C$ are finite dimensional $\ku$-linear bundles over
the vertical groupoid $\B \rightrightarrows \Hc$; that is,
$\Hc$-graded vector spaces endowed with a left action of the
vertical groupoid $\B \rightrightarrows \Hc$ by linear
isomorphisms, \emph{cf.} \cite{AN}. Here, a box $A \in \B$ acts on
$V = \oplus_{x \in \Hc}V_x$ via the linear isomorphism $A:
V_{b(A)} \to V_{t(A)}$.

\bigbreak
There is a structure of rigid monoidal category on $\Rep \ku \T$.
The unit object is the target subalgebra $\ku \T_t = \ku(\uno_{E}:
E \in \Ec) \simeq (\ku \Ec)^{\op}$. If $V$, $W$ are $\T$-bundles
then $V \otimes W := V \otimes_{\Ec^{\op}}W$, where $\ku\Ec^{\op}$
acts on the left via the target subalgebra and on the right via
the left action of the source subalgebra. The action of $\B$ on $V
\otimes W$ is given by $\Delta$.  The \emph{dual} $V^*$ of an
object $V = \oplus_{x \in \Hc}V_x \in \C$ has $\Hc$-grading
$(V^*)_x = (V_{x^{-1}})^*$, $x \in \Hc$. If  $A \in \B$, then it
acts by the transpose of $A^{-1} : (V^*)_{b(A)} \to (V^*)_{t(A)}$.

\bigbreak The core groupoids $\Dc$ and $\Ec$ are embedded ``diagonally" into the algebra $\ku\T$. The multiplication between elements of $\Dc$ or $\Ec$ and those in $\B$ is related to the actions in Proposition \ref{acciones-core}.

\begin{lema} Let $A \in \B$, $D \in \Dc$, $E \in \Ec$. Then

\begin{align}\pagebreak\label{D1A}\pagebreak {}_D\uno. A &=
\begin{cases} D \rightharpoonup A, \quad \text{if } rt(A) = e(D),
\\ 0, \quad \text{otherwise;}
\end{cases}
\\ \label{AD1} \pagebreak A.{}_D\uno &=
\begin{cases} A \leftharpoonup D, \quad \text{if } rb(A) = s(D),
\\ 0, \quad \text{otherwise;}
\end{cases}\end{align}
\begin{align} \label{1EA}\pagebreak \uno_E . A &= \begin{cases} A
\leftharpoondown E, \quad \text{if } lt(A) = s(E),
\\ 0, \quad \text{otherwise;}
\end{cases}
\\ \label{A1E} \pagebreak A. \uno_E &= \begin{cases} E \rightharpoondown A,
\quad \text{if } lb(A) = e(E),
\\ 0, \quad \text{otherwise.}
\end{cases}\end{align} \qed \end{lema}

 Let $D\in \Dc$, $E\in \Ec$, $z,w\in \Hc$. Then, by the Lemma above, we have
\begin{align}\label{D1x} {}_D\uno.  \iddh \, z
&= \begin{cases} \{\iddh\, z D\} = \iddh \, (zt(D))
. {}_D\uno, \quad \text{if } r(z) = e(D),
\\ 0, \quad \text{otherwise;}
\end{cases}
\\ \label{1Ex} \uno_E. \iddh \, w
&= \begin{cases} \{E\iddh\, b(E)^{-1}w\} = \iddh \, (b(E)^{-1}w) . {}_E\uno, \quad \text{if } l(w) = e(E),
\\ 0, \quad \text{otherwise.}
\end{cases}\end{align}

If $D \in \Dc$, set $\theta(D) := \theta(e(D)) = \theta(s(D))$, \emph{cf.} Lemma \ref{d-cnt}. If $E \in \Ec$, set $\theta(E) := \theta(E^{-1})$.

\begin{prop}
\begin{align}
\label{delta-1} \Delta(1) &= \sum_{D\in \Dc}\frac1{\theta(D)} \,
{}_D\uno \otimes \uno_{D^{\daga}}
\\ \label{delta-1bis} &= \sum_{E\in \Ec}\frac1{\theta(E)} \,
{}_{E^{\daga}}\uno \otimes \uno_E;
\\
 \label{delta-x} \Delta(\iddh \, x) &= \Delta(1)\sum_{zw = x}
\iddh \, z\otimes \iddh \, w, \qquad x\in \Hc.
\end{align}
\end{prop}
\pf Let $x\in \Hc$. We compute
\begin{align*}
&\Delta(\iddh \, x) = \sum_{AB = \iddh \, x}
\frac1{\urcorner(B)} \; A\otimes B \\ &=
\sum_{D\in \Dc, \, z \in \Hc: r(z) = e(D)}\;
\sum_{y \in \Hc, \, l(y) = e(E); \{\iddv z \, DE \, \iddv y\} = \iddv x}
\frac1{\theta(D)}\{\iddv z \, D\} \otimes \{E \, \iddv y\}
\\ &=
\sum_{D\in \Dc, \, z \in \Hc: r(z) = e(D)} \frac1{\theta(D)}\{\iddv z \, D\} \otimes \{D^{\daga} \iddv (t(D)^{-1}z^{-1}x)\}.
\end{align*}
Here the first equality is by definition; the second uses the bijections \eqref{b-left} and \eqref{b-right}; the third, Remark
\ref{daga-appl}. This implies \eqref{delta-1}, and \eqref{delta-1bis} follows by a routine change of variables. Starting from the last equation and using \eqref{D1x}, we have
\begin{align*}
\Delta(\iddh \, x) &=
\sum_{D\in \Dc, \, z \in \Hc: r(z) = e(D)} \frac1{\theta(D)}{}_D\uno.  \iddh \, z \otimes \uno_{D^{\daga}}. \iddv (z^{-1}x),
\end{align*}
since $t(D) = b(D^{\daga})$. This implies \eqref{delta-x}.  \epf

As a consequence, if $V$, $W$ are $\T$-bundles then the homogeneous components of the tensor product are given by $$(V \otimes W)_x = \Delta(1) \big(\sum_{zw = x} V_z\otimes W_w\big),\qquad x\in \Hc.$$

\section{Cocycle deformations}\label{extcocycles}

\subsection{Generalities}

\

Let $\T$ be a  \emph{finite}
double groupoid. As before,  let $\ku\T$ denote the  $\ku$-vector
space with basis $\B$. We shall consider the following structures
on $\ku\T$.

\bigbreak\subsection*{Algebra structure}

\

We  deform the groupoid
algebra structure on $\ku\T$ corresponding to the vertical composition
groupoid $\B\rightrightarrows \Hc$.

\begin{lema} Let $\sigma: \B \Times{b}{t} \B \to \ku^{\times}$
be a function and define a multiplication in $\ku\T$ by
\begin{equation}
\label{producto} A.B = \begin{cases} \sigma(A, B)\,
\begin{matrix}A \vspace{-7pt}\\B\end{matrix}, \quad
\text{if } \displaystyle\frac{A}{B}, \\
0 , \qquad\qquad\quad \text{otherwise}. \end{cases} \qquad
\text{for all }A, B \in \B.
\end{equation}
This multiplication is associative if and only if $\sigma$ is a vertical 2-cocycle:
\begin{equation}
\label{cociclo-sigma} \sigma(A, B) \sigma\left(\begin{matrix}A
\vspace{-7pt}\\B\end{matrix}, C\right) = \sigma(B, C)
\sigma\left(A,\begin{matrix}B \vspace{-7pt}\\C\end{matrix}\right),
\end{equation}
for all $A, B, C \in \B$: $\dfrac{\displaystyle\dfrac{A}{B}}{C}$.
If this happens,  there is a unit $$\unosigma := \sum_{x \in \Hc}
\displaystyle\frac{1}{\sigma(\iddv \,{x}, \iddv \,{x})} \,\iddv
\,{x}.$$
\end{lema}
\pf The proof of the first claim is routine. For the second, the
identities $\sigma(A, \iddv b(A)) = \sigma(\iddv b(A), \iddv
b(A))$, $\sigma(\iddv t(A), A) = \sigma(\iddv t(A), \iddv t(A))$
are needed, but these follow from \eqref{cociclo-sigma}. \epf

If \eqref{cociclo-sigma} holds, the unit is $\uno := \sum_{x
\in \Hc} \iddv \,{x}$ if and only if $\sigma$ is normalized:
\begin{equation}\label{norm-sigma}
 \sigma(A, \id \, {t(A)}) = \sigma(\id \, {b(A)}, A) = 1,
\text{ for all }A\in \B.\end{equation}
Up to a change of basis, one can always assume that $\sigma$ is normalized.

\bigbreak
\subsection*{Coalgebra structure}

\

Dually, we shall deform the
groupoid coalgebra structure on $\ku\T$ corresponding to the
horizontal  composition groupoid $\B \rightrightarrows \Vc$.

\begin{lema} Let $\tau: \B \Times{r}{l} \B \to \ku^{\times}$
be a function and define a comultiplication in $\ku\T$ by
\begin{equation}
\label{coproducto} \Delta(A) = \sum_{A = B C} \tau(B, C) \,
B\otimes C, \qquad A \in \B.
\end{equation}
This comultiplication is coassociative if and only if
\begin{equation}
\label{cociclo-tau} \tau(A, B) \tau(AB, C) = \tau(B, C)
\tau(A,BC),
\end{equation}
for all $A, B, C \in \B$: $A\vert B\vert C$. If this happens, then
there is a counit $\varepsilon^{\tau}: \ku\T \to \ku$ given by
$$\varepsilon^{\tau}(A) =
\begin{cases} \displaystyle\frac{1}{\tau(A, A)},
\quad \text{if } A \in \Vc,  \\ 0 , \quad \text{otherwise}.
\end{cases}$$ \qed\end{lema}

\bigbreak
\subsection*{Weak Hopf algebra structure}

\

We fix
$\sigma: \B \Times{b}{t} \B \to \ku^{\times}$ satisfying
\eqref{cociclo-sigma} and \eqref{norm-sigma}, and $\tau: \B \Times{r}{l} \B \to
\ku^{\times}$ satisfying \eqref{cociclo-tau}.
We denote by $\ku^{\tau}_{\sigma}\T$  the vector space $\ku\T$ with multiplication \eqref{producto} and  comultiplication
\eqref{coproducto}. Observe that it is not possible to normalize $\sigma$ and $\tau$ simultaneously. We begin by the following straightforward proposition.

\begin{prop}\label{general}   $\ku^{\tau}_{\sigma}\T$ is a weak bialgebra
if and only if the following hold:

\begin{equation}
\label{multiplicativa} \sigma(A, B)\tau(X,Y) = \sum \sigma(U,
R)\sigma(V, S)\tau(U,V)\tau(R,S),
\end{equation}
for all $A, B, X,
Y \in \B: XY = \begin{matrix}A\vspace{-4pt}\\B\end{matrix}$, where
the index set is $[X, Y, A, B]$;
\begin{align} \label{2.2I} \tau\left(A, \begin{tiny}\begin{matrix}U\vspace{-4pt}\\V\end{matrix}
\end{tiny}\right)
\tau\left(A\, \begin{tiny}\begin{matrix}U\vspace{-4pt}\\V\end{matrix}\end{tiny},
C\right) &= \tau(A, U)
 \tau(V,C) \sigma(U,V),
\\ \label{2.2II} \tau\left(A, \begin{tiny}
\begin{matrix}W\vspace{-4pt}\\Z\end{matrix}\end{tiny}\right)
\tau\left(A\,\begin{tiny}\begin{matrix}W\vspace{-4pt}\\Z\end{matrix}\end{tiny}, C\right) &= \tau(A, Z)
\tau(W,C) \sigma(W,Z),\end{align}
for all $A, C, U,
V, W, Z \in \B:$ $AU, VC, AZ, WC \in \Hc$;
\begin{align} \label{2.3I} \sigma(A,UV)
\sigma\left(\begin{tiny}\begin{matrix}A\vspace{-4pt}\\UV\end{matrix}\end{tiny},C\right)
\tau\left(\begin{tiny}\begin{matrix}A\vspace{-4pt}\\U\end{matrix}\end{tiny},
\begin{tiny}\begin{matrix}A\vspace{-4pt}\\U\end{matrix}\end{tiny}\right)
\tau\left(\begin{tiny}\begin{matrix}V\vspace{-4pt}\\C\end{matrix}\end{tiny},
\begin{tiny}\begin{matrix}V\vspace{-4pt}\\C\end{matrix}\end{tiny}\right)
&= \tau(U,V)
\sigma(A,U) \sigma(V,C)
\tau\left(\begin{tiny}\begin{matrix}A\vspace{-4pt}\\UV\vspace{-4pt}\\C\end{matrix}\end{tiny},
\begin{tiny}\begin{matrix}A\vspace{-4pt}\\UV\vspace{-4pt}\\C\end{matrix}\end{tiny}\right),
\end{align}
\begin{multline} \label{2.3II} \sigma(A,WZ)
\sigma\left(\begin{tiny}\begin{matrix}A\vspace{-4pt}\\WZ\end{matrix}\end{tiny},C\right)
\tau\left(\begin{tiny}\begin{matrix}A\vspace{-4pt}\\Z\end{matrix}\end{tiny},
\begin{tiny}\begin{matrix}A\vspace{-4pt}\\Z\end{matrix}\end{tiny}\right)
\tau\left(\begin{tiny}\begin{matrix}W\vspace{-4pt}\\C\end{matrix}\end{tiny},
\begin{tiny}\begin{matrix}W\vspace{-4pt}\\C\end{matrix}\end{tiny}\right)
\\= \tau(W,Z)
\sigma(A,Z) \sigma(W,C)
\tau\left(\begin{tiny}\begin{matrix}A\vspace{-4pt}\\WZ\vspace{-4pt}\\C\end{matrix}\end{tiny},
\begin{tiny}\begin{matrix}A\vspace{-4pt}\\WZ\vspace{-4pt}\\C\end{matrix}\end{tiny}\right),
\end{multline}
for all $A, C, U,
V, W, Z \in \B:$ $\begin{tiny}\begin{matrix}A\vspace{-4pt}\\U\end{matrix}\end{tiny}$,
$\begin{tiny}\begin{matrix}V\vspace{-4pt}\\C\end{matrix}\end{tiny}$,
$\begin{tiny}\begin{matrix}A\vspace{-4pt}\\Z\end{matrix}\end{tiny}$,
$\begin{tiny}\begin{matrix}W\vspace{-4pt}\\C\end{matrix}\end{tiny}
\in \Vc$.

In this case, the source and target maps are given, respectively, by
\begin{align} \label{sourcegen} \epsilon_s(A) & =
\begin{cases}
\sum_{x \in \Hc, \, r(x) = br(A)}
\dfrac{\tau(\iddv x \, \phi(A),  \phi(A)^h) \sigma(A,\phi(A)^h) }
{\tau(\iddv r(A), \iddv r(A))}
\{\iddv x \, \phi(A)\}, \\ \qquad\qquad \qquad\qquad\qquad\qquad\qquad\qquad\qquad\qquad
\text{if} \quad t(A) \in \Pc, \\
0, \qquad \text{otherwise};  \end{cases}\end{align}  \begin{align}
\label{targetgen} \epsilon_t(A) & =
\begin{cases} \sum_{y \in \Hc, \, l(y) = bl(A)}
\dfrac{\tau(\psi(A)^h,  \psi(A)\iddv y) \sigma(\psi(A)^h, A) }
{\tau(\iddv l(A), \iddv l(A))}
\{\psi(A) \iddv y \}, \\ \qquad\qquad \qquad\qquad\qquad\qquad\qquad\qquad\qquad\qquad
\text{if} \quad b(A)  \in \Pc, \\
0, \qquad \text{otherwise}.  \end{cases}  \end{align}
\end{prop}

\pf  We first show:

\begin{step}\label{multiplicativa-l}
The comultiplication \eqref{coproducto} is
multiplicative with respect to the multiplication \eqref{producto}
if and only if
\eqref{multiplicativa} holds.
\end{step}

Let $A, B
 \in \B$. It follows from the definitions that $$\Delta(A.B) =
\sum_{\begin{tiny}XY =
\begin{matrix}A \vspace{-4pt}\\B\end{matrix}\end{tiny}}
\sigma(A, B)\tau(X,Y) \quad  X \otimes Y.$$ On the other hand,
$$\Delta(A) . \Delta(B) = \sum_{}\sigma(U,
R)\sigma(V, S)\tau(U,V)\tau(R,S) \,  \begin{matrix}U \vspace{-4pt}\\R\end{matrix}
\otimes \begin{matrix}V \vspace{-4pt}\\S\end{matrix},$$ where the
sum runs over all elements $U, V, R, S \in \B$, such that
\begin{tiny}\begin{tabular}{p{0,2cm}|p{0,2cm}} $U$ & $V$
\\ \hline $R$ & $S$ \end{tabular}\end{tiny}, $UV = A$ and  $RS = B$.
It is thus clear that $\Delta(A) . \Delta(B) = 0 = \Delta(A . B)$,
if $A$ and $B$ are not vertically composable. If
$\displaystyle\frac{A}{B}$, then \eqref{d-mult} is
equivalent to \eqref{multiplicativa}.

\medbreak
\begin{step}\label{ax-unit-l}
The first, resp. the second,
equality in  \eqref{ax-unit} is
equivalent to
\eqref{2.2I}, resp. \eqref{2.2II}.
\end{step}

We
have \begin{align*}\Delta^{(2)} (\uno )  = \sum_{ABC \in \Hc}
\tau\left(A, B\right)
\tau\left(AB, C\right) \, A \otimes B \otimes C.\end{align*}
On the other hand, by \eqref{norm-sigma}, we have
\begin{align*} \left( \Delta(\uno) \otimes \uno \right)
&\left( \uno \otimes \Delta(\uno) \right) \\ &=
\sum_{\begin{tiny}\begin{matrix}U\\ V \end{matrix}\end{tiny}
= B, \, AU \in \Hc, \, VC \in \Hc}
\tau(A, U)
\tau(V,C) \sigma(U,V) \quad A \otimes
\begin{matrix} U \vspace{-4pt}\\ V \end{matrix} \otimes C \\
\left( \uno \otimes \Delta(\uno) \right)
&\left( \Delta(\uno) \otimes \uno \right) \\ &=
\sum_{\begin{tiny}\begin{matrix}W\\ Z
\end{matrix}\end{tiny} = B, \, AZ \in \Hc, \, WC \in \Hc}
\tau(A, Z)
\tau(W,C) \sigma(W,Z) \quad A \otimes
\begin{matrix} W \vspace{-4pt}\\ Z \end{matrix} \otimes C.\end{align*}
Thanks to the equivalences (i) $\Longleftrightarrow$ (ii) and (i)
$\Longleftrightarrow$ (iii) in Lemma \ref{l-unit}, the claim in
Step 2 follows.

\medbreak
\begin{step}\label{ax-counit-l}
The first, resp. the second,
equality in  \eqref{ax-counit} is
equivalent to
\eqref{2.3I}, resp. \eqref{2.3II}.
\end{step}
Let  $A, B, C \in \B$. We have
\begin{equation*}\epsilon(A.B.C) = \begin{cases} \sigma(A,B)
\sigma\left(\begin{tiny}\begin{matrix}A\vspace{-4pt}\\B\end{matrix}\end{tiny},C\right)
\tau\left(\begin{tiny}\begin{matrix}A\vspace{-4pt}\\B\vspace{-4pt}\\C\end{matrix}\end{tiny},
\begin{tiny}\begin{matrix}A\vspace{-4pt}\\B\vspace{-4pt}\\C\end{matrix}\end{tiny}\right)^{-1}, \quad
\text{if } \quad
\begin{tiny}\begin{matrix}A \vspace{-4pt}\\ B \vspace{-4pt}\\ C \end{matrix}\end{tiny} \; \in \Vc,  \\ 0 , \quad \text{otherwise}.
\end{cases} \end{equation*}
On the other hand, we have

\begin{align*}\epsilon(A.B_1)&\epsilon(B_2.C)
= \sum_{UV = B} \tau(U, V) \, \epsilon(A.U)\epsilon(V.C)
\\  &=  \sum_{UV = B, \,
\begin{tiny}\begin{matrix} A\\ U \end{matrix}\end{tiny} \in \Vc, \,
\begin{tiny}\begin{matrix} V\\ C \end{matrix}\end{tiny} \in \Vc}
\tau(U, V)\sigma(A,U)\sigma(V,C)\tau\left(\begin{tiny}\begin{matrix}A\\U\end{matrix}\end{tiny},
\begin{tiny}\begin{matrix}A\\U\end{matrix}\end{tiny}\right)^{-1}
\tau\left(\begin{tiny}\begin{matrix}V\\C\end{matrix}\end{tiny},
\begin{tiny}\begin{matrix}V\\C\end{matrix}\end{tiny}\right)^{-1}. \end{align*}
By Lemma \ref{l-counit}, there exist
unique such $U$ and $V$  in the index set of the last sum.
Thus the first equality in  \eqref{ax-counit} is
equivalent to
\eqref{2.3I}. Similarly, the second
equality in  \eqref{ax-counit} is
equivalent to
\eqref{2.3II}.

\medbreak
\begin{step}\label{sourcetarg-l}
The  source and target maps.
\end{step}

\medbreak  Let $A \in \B$. We compute
\begin{equation}\label{ep-s}\epsilon_s(A)  = \sum_{XY \in \Hc}
\tau(X,Y) \; X \; \epsilon(A.Y) =
\sum_{XY \in \Hc, \,
\begin{tiny}\begin{matrix} A\vspace{-4pt} \\Y\end{matrix}\end{tiny} \in \Vc} \;
\frac{\tau(X,Y)\sigma(A,Y)}
{\tau(\begin{tiny}\begin{matrix} A\vspace{-4pt} \\Y\end{matrix}\end{tiny},
\begin{tiny}\begin{matrix} A\vspace{-4pt} \\Y\end{matrix}\end{tiny})}  \quad X.
\end{equation}
By Lemma \ref{conteo-st} (i), $Y = \phi(A)^h$, $X = \iddv x
\phi(A)$, for a unique $x \in \Hc$ such that $r(x) = rb(A)$,
whence formula \eqref{sourcegen}.

As for the target map, we have
\begin{equation}\label{ep-t}\epsilon_t(A)
= \sum_{XY \in \Hc} \tau(X,Y) \; \epsilon(X.A) Y =
\sum_{XY \in \Hc, \,
\begin{tiny}\begin{matrix}X \vspace{-4pt}\\ A \end{matrix}\end{tiny}
\in \Vc}
\frac{\tau(X,Y)\sigma(X,A)}
{\tau(\begin{tiny}\begin{matrix} X\vspace{-4pt} \\A\end{matrix}\end{tiny},
\begin{tiny}\begin{matrix} X\vspace{-4pt} \\A\end{matrix}\end{tiny})}  \; Y.
\end{equation}
By Lemma
\ref{conteo-st} (ii), \eqref{targetgen} follows.
\epf

We denote by $\ku^{\tau}\T$  the algebra and coalgebra
$\ku^{\tau}_{\sigma}\T$ when $\sigma = 1$ is the trivial cocycle.

\begin{cor}\label{sigmatrivial}   $\ku^{\tau}\T$ is a weak bialgebra
if and only if

\begin{equation}
\label{multiplicativa-st} \tau(X,Y) = \sum \tau(U,V)\tau(R,S),
\end{equation}
for all $A, B, X,
Y \in \B: XY = \begin{matrix}A\vspace{-4pt}\\B\end{matrix}$, where
the index set is $[X, Y, A, B]$;
\begin{align} \label{2.2I-st}
\tau\left(A, \begin{tiny}\begin{matrix}U\vspace{-4pt}\\V\end{matrix}
\end{tiny}\right)
\tau\left(A\, \begin{tiny}\begin{matrix}U\vspace{-4pt}\\V\end{matrix}\end{tiny},
C\right) &= \tau(A, U)
 \tau(V,C) ,
\\ \label{2.2II-st} \tau\left(A, \begin{tiny}
\begin{matrix}W\vspace{-4pt}\\Z\end{matrix}\end{tiny}\right)
\tau\left(A\,\begin{tiny}\begin{matrix}W\vspace{-4pt}\\Z\end{matrix}\end{tiny}, C\right)
&= \tau(A, Z)
\tau(W,C),
\end{align}
for all $A, C, U,
V, W, Z \in \B:$ $AU, VC, AZ, WC \in \Hc$;
\begin{align}
\label{2.3Ist}
\tau\left(\begin{tiny}\begin{matrix}A\vspace{-4pt}\\U\end{matrix}\end{tiny},
\begin{tiny}\begin{matrix}A\vspace{-4pt}\\U\end{matrix}\end{tiny}\right)
\tau\left(\begin{tiny}\begin{matrix}V\vspace{-4pt}\\C\end{matrix}\end{tiny},
\begin{tiny}\begin{matrix}V\vspace{-4pt}\\C\end{matrix}\end{tiny}\right)
&= \tau(U,V)
\tau\left(\begin{tiny}\begin{matrix}A\vspace{-4pt}\\B\vspace{-4pt}\\C\end{matrix}\end{tiny},
\begin{tiny}\begin{matrix}A\vspace{-4pt}\\B\vspace{-4pt}\\C\end{matrix}\end{tiny}\right),
\\ \label{2.3II-st}
\tau\left(\begin{tiny}\begin{matrix}A\vspace{-4pt}\\Z\end{matrix}\end{tiny},
\begin{tiny}\begin{matrix}A\vspace{-4pt}\\Z\end{matrix}\end{tiny}\right)
\tau\left(\begin{tiny}\begin{matrix}W\vspace{-4pt}\\C\end{matrix}\end{tiny},
\begin{tiny}\begin{matrix}W\vspace{-4pt}\\C\end{matrix}\end{tiny}\right)
&= \tau(W,Z)
\tau\left(\begin{tiny}\begin{matrix}A\vspace{-4pt}\\B\vspace{-4pt}\\C\end{matrix}\end{tiny},
\begin{tiny}\begin{matrix}A\vspace{-4pt}\\B\vspace{-4pt}\\C\end{matrix}\end{tiny}\right),
\end{align}
for all $A, C, U, V, W, Z \in \B:$ $\mvert{A}{U}$,
$\begin{tiny}\begin{matrix}V\vspace{-4pt}\\C\end{matrix}\end{tiny}$,
$\begin{tiny}\begin{matrix}A\vspace{-4pt}\\Z\end{matrix}\end{tiny}$,
$\begin{tiny}\begin{matrix}W\vspace{-4pt}\\C\end{matrix}\end{tiny}
\in \Vc$. \qed
\end{cor}

\bigbreak
\subsection{Generalized corner functions}\label{gcf}

\

We consider now the following particular case.
We fix a function $\vartheta: \Pc \to \ku^{\times}$ and introduce
$\cx: \B \to \ku^{\times}$, $\tau: \B \Times{r}{l} \B \to \ku^{\times}$ by
$$
\cx(B) = \vartheta(bl(B)), \qquad \tau(A, B) = \cx(B).
$$
It is immediate that
\begin{align}\label{basico-cxh}\cx(XY) &= \cx(X),
\\ \label{basico-cxv}
\cx(\mvert{U}{V}) &= \cx(V),\end{align}
for all $X,Y,U,V  \in \B$ appropriately composable. Note that condition \eqref{cociclo-tau} for $\tau$ is equivalent to condition \eqref{basico-cxh} for $\cx$. Also, condition \eqref{2.2II-st} for $\tau$ is equivalent to condition \eqref{basico-cxv} for $\cx$.

\begin{obs} Let $\cx: \B \to \ku^{\times}$ be any function and define $\tau: \B \Times{r}{l} \B \to \ku^{\times}$ by $\tau(A, B) = \cx(B)$. Suppose that $\cx$ satisfies \eqref{basico-cxh}  and \eqref{basico-cxv}. Then $\cx(B) = \vartheta(bl(B))$ for any $B\in \B$, where $\vartheta: \Pc \to \ku^{\times}$ is given by $\vartheta(P) = \cx(\Theta_P)$, $P\in \Pc$.
\end{obs}

To stress the dependence on $\vartheta$, we shall use the notation
$\ku^{(\vartheta)}\T$ instead of $\ku^{\tau}\T$.

\begin{theorem}\label{cx}   $\ku^{(\vartheta)}\T$ is a weak bialgebra
if and only if

\begin{equation}
\label{multiplicativa-cx} 1 = \sum_{V\in \B: t(V) = x, \, r(V) = g} \cx(V),
\end{equation}
for all $(g, x) \in \Vc \Times{t}{r}\Hc$.

\medbreak When this condition holds, $\ku^{(\vartheta)}\T$ is a
semisimple weak Hopf algebra, with antipode determined by the formula

\begin{equation}\label{antipodegen-cx}\mathcal S(A) =
\dfrac{\cx(A^{-1})}{\cx(A^h)}\; A^{-1}
=\dfrac{\vartheta(tr(A))}{\vartheta(br(A))}\; A^{-1},\end{equation}

\noindent for all $A
\in \B$.  The source and target maps are given, respectively, by
\begin{align} \label{sourcegen-cx} \epsilon_s(A) & =
\begin{cases}
{}_{\phi(A)}\uno, \qquad
\text{if} \quad t(A) \in \Pc, \\
0, \qquad \text{otherwise};  \end{cases}  \\
\label{targetgen-cx} \epsilon_t(A) & =
\begin{cases} \dfrac{\cx(A^{-1})}{\cx(A^h)}\; \uno_{\psi(A)}\qquad
\text{if} \quad b(A)  \in \Pc, \\
0, \qquad \text{otherwise}.  \end{cases}
\end{align}
The source and target subalgebras are isomorphic to the groupoid
algebras $\ku \Dc$ and  $(\ku \Ec)^{\op}$, respectively.

\end{theorem}

Note that \eqref{multiplicativa-cx} implies that $\T$ satisfies the filling condition \eqref{neq0}.

\pf Clearly, \eqref{multiplicativa-st} is equivalent to
\begin{equation*}
 \cx(Y) = \tau(X,Y) = \sum \tau(U, V)\tau(R, S)
= \sum \cx(V)\cx(S)
\end{equation*}
for all $A, B, X,
Y \in \B: XY = \begin{matrix}A\vspace{-4pt}\\B\end{matrix}$, where
the index set is $[X, Y, A, B]$. Since $\cx(Y) = \cx(S)$, we see
that \eqref{multiplicativa-cx} is equivalent to
\eqref{multiplicativa-st}.

\medbreak
Now conditions \eqref{2.2II-st}  and \eqref{2.3II-st} follow at once from
\eqref{basico-cxv}, while \eqref{2.2I-st} and \eqref{2.3Ist} follow easily using
that $VC \in \Hc$, resp. that $\mvert{A}{U} \in \Vc$.

\medbreak
It remains to consider the axioms of the antipode.
We first show \eqref{atp-1}. Let $A \in \B$. Then
\begin{equation*}
m(\id \otimes \Ss) \Delta (A)  =
\sum_{\frac{U}{V^{-1}}, \, UV = A} \tau(U,V)\dfrac{\cx(V^{-1})}{\cx(V^h)} \;
\begin{matrix}U \vspace{-4pt}\\ \, V^{-1}\end{matrix}; \end{equation*}
thus this sum vanishes, and so equation \eqref{atp-1} is true,
unless $b(A) \in \Pc$. We shall now assume that this is the case.
Letting $Y = \begin{matrix}U
\vspace{-4pt}\\ \, V^{-1}\end{matrix}$, and $X =
\left\{\begin{matrix}\iddv l(A) \vspace{-4pt}\\ A^v
\end{matrix}\right\}$, the last sum equals
\begin{align*}
&\sum_{XY \in \Hc, \; \begin{tiny}\begin{matrix}X \\
A \end{matrix}\end{tiny} \in \Vc} \quad \sum_{UV = A, \,
\begin{tiny}\begin{matrix}U \\ \, V^{-1}\end{matrix}\end{tiny} =
Y} \cx(V) \dfrac{\cx(V^{-1})}{\cx(V^h)}\, Y
\\ &=
 \sum_{XY \in \Hc,
\,\begin{tiny}
\begin{matrix}X \\ A \end{matrix}\end{tiny} \in \Vc}
\quad \dfrac{\cx(Y)}{\cx(A)}\sum_{UV = A, \,
\begin{tiny}\begin{matrix}U \\ \, V^{-1}\end{matrix}\end{tiny} =
Y} \cx(V) \, Y
\\ &=
 \sum_{XY \in \Hc,
\,\begin{tiny}
\begin{matrix}X \\ A \end{matrix}\end{tiny} \in \Vc}
\quad \dfrac{\cx(Y)}{\cx(A)}
 Y \\ & =
\sum_{XY \in \Hc, \,
\begin{tiny}\begin{matrix}X \vspace{-4pt}\\ A \end{matrix}\end{tiny}
\in \Vc}
\frac{\tau(X,Y)}
{\tau(\begin{tiny}\begin{matrix} X\vspace{-4pt} \\A\end{matrix}\end{tiny},
\begin{tiny}\begin{matrix} X\vspace{-4pt} \\A\end{matrix}\end{tiny})}  \; Y,\end{align*}
and we are done by \eqref{ep-t}.
We have used that $\cx(Y) = \cx(V^{-1})$
and $\cx(V^h)= \vartheta(br(V))= \vartheta(tl(V^{-1})) = \vartheta(bl(U)) = \cx(A).
$

Relation \eqref{atp-2} is proved similarly: let $A \in \B$. Then
\begin{equation*}m(\Ss \otimes \id) \Delta (A)  =
\sum_{\frac{U^{-1}}{V}, \; UV = A} \cx(V)
\dfrac{\cx(U^{-1})}{\cx(U^h)} \quad \begin{matrix}\quad
U^{-1} \vspace{-4pt}\\ \, V\end{matrix}; \end{equation*} as
before, this sum vanishes, implying equation \eqref{atp-2}, in the
case where $t(A)$ is not in $\Pc$. We shall assume that $t(A) \in
\Pc$. Letting $X = \begin{matrix}U^{-1} \vspace{-4pt}\\  V\quad
\end{matrix}$ and $Y = \left\{\begin{matrix}A^v \vspace{-4pt}\\
\iddv \, r(A)\end{matrix}\right\}$, we get
\begin{align*} m(\Ss \otimes \id) \Delta (A)  &=\sum_{XY \in \Hc, \,\begin{tiny}
\begin{matrix}A \\ Y \end{matrix}\end{tiny} \in \Vc}
\quad \sum_{UV = A,\;
\begin{tiny}\begin{matrix}U^{-1} \\ V \quad\end{matrix}\end{tiny}
= X} \cx(U^{-1})\quad X,  \end{align*} which gives
\eqref{atp-2} in view of \eqref{ep-s}, after observing that
$\cx(U^h) = \cx(V)$.

We finally show \eqref{atp-3}. Let $A \in \B$. We have
\begin{align*}m^{(2)}(\Ss \otimes \id \otimes \Ss) \Delta^{(2)}(A)
& = \sum_{XYZ = A}
\dfrac{\cx(Y)\cx(Z)\cx(X^{-1})\cx(Z^{-1})}{\cx(X^{h})\cx(Z^{h})}
 \quad X^{-1}.Y.Z^{-1} \\ & =
\Big(\sum_{} \;
\dfrac{\cx(Y)\cx(Z)\cx(X^{-1})}{\cx(X^{h})}
\Big)
\dfrac{\cx(A^{-1})}{\cx(A^{h})}
 \quad A^{-1},
\end{align*} where the last sum runs over all triples $(X, Y,
Z)$ satisfying \eqref{tricot1}, \eqref{tricot2}; see Lemma
\ref{l-atp3}. Here we have used that $\cx(A^{h}) = \cx(Z^{h})$, which follows from
\eqref{tricot1}, and $\cx(A^{-1}) = \cx(Z^{-1})$, which follows from
\eqref{tricot2}. Now
$\cx(Y) = \vartheta(bl(Y)) = \vartheta(br(X)) = \vartheta(bl(X^h)) = \cx(X^{h})$
and thus, because of Lemma \ref{conteo-ant} and hypothesis \eqref{multiplicativa-cx}, we get
\begin{multline*}\sum_{} \;
\dfrac{\cx(Y)\cx(Z)\cx(X^{-1})}{\cx(X^{h})}
= \Big(\sum_{Z\in \B: t(Z) = t(A), \, r(Z) = r(A)} \; \cx(Z)\Big) \\ \times
\Big(\sum_{X\in \B: b(X) = b(A), \, l(X) = l(A)} \; \cx(X^{-1})\Big) = 1,
\end{multline*} proving \eqref{atp-3}.

\medbreak Finally, the maps $\ku \Dc \to (\ku^{(\vartheta)}
\T)_s$, $D \mapsto {}_D\uno$, and $(\ku \Ec)^{\op} \to
(\ku^{(\vartheta)} \T)_t$, $E \mapsto \uno_E$, determine algebra
isomorphisms, by Lemma \ref{gpd-st}. \epf

\bigbreak
\begin{obs} Let $\vartheta: \Pc \to \ku^{\times}$ be any function and define $\cx$ as before. For any $(g, x) \in \Vc \Times{t}{r}\Hc$, set $$c(g,x):= \sum_{V\in \B: t(V) = x, \, r(V) = g} \cx(V).$$
Let us say that $P\in \Pc$ is \emph{affiliated} to $(g, x)$ if there exists $V\in \B$: $t(V) = x$, $r(V) = g$, $bl(V) = P$.
Assume that the filling condition \eqref{neq0} holds and that
\begin{equation}\label{normalization-cx}
c(g,x) = c(\id_{\Vc} P, \id_{\Hc} P) \neq 0
\end{equation}
for all $(g, x)$ such that $P$ is affiliated to $(g, x)$.
Let $\widetilde\vartheta: \Pc \to \ku^{\times}$ be given by $$\widetilde\vartheta(P) = \dfrac{1}{c(\id_{\Vc} P, \id_{\Hc} P)}\vartheta(P), \qquad p\in \Pc.$$
Then the $\widetilde\cx$ corresponding to $\widetilde\vartheta$ satisfies \eqref{multiplicativa-cx}.
\end{obs}

\bigbreak
{\it Proof of Theorem \ref{bicross}.} Consider $\vartheta: \Pc\to \ku^{\times}$,
$\vartheta(P) = \dfrac{1}{\theta(P)}$. We only need to show that
$\vartheta$ satisfies \eqref{multiplicativa-cx}.
Let $(g, x) \in \Vc \Times{t}{r}\Hc$. Because of the filling condition
\eqref{neq0}, there exists $A\in \B$ such that $t(A) = x$, $r(A) = g$;
Set $Y=A$, $X= \iddv\, l(A)$,$X= \iddh\, b(A)$. Then,
by Proposition \ref{rel-corners}, we have
\begin{align*}
\sum_{V\in \B: t(V) = x, \, r(V) = g} \; \cx(V) &=
 \sum_{UV = A, \, RS = B,
\begin{tiny}\begin{matrix}U\\R\end{matrix}\end{tiny} = X,
\, \begin{tiny}\begin{matrix}V\\ S\end{matrix}\end{tiny} =
Y}\dfrac{1}{\urcorner(V)} \\
&= \sum_{UV = A, \, RS = B,
\begin{tiny}\begin{matrix}U\\R\end{matrix}\end{tiny} = X, \,
\begin{tiny}\begin{matrix}V \\ S\end{matrix}\end{tiny} =
Y}\dfrac{1}{\urcorner(t(Y), r(A))} = 1.
\end{align*}
\qed

\bigbreak
Let $\tau$ be given by a general $\vartheta$ satisfying
\eqref{multiplicativa-cx}. The following statement generalizes
Proposition \ref{regular}.

\begin{prop}\label{regular-gen} $\ku^{(\vartheta)} \T$ is a regular weak Hopf algebra
if and only if $\vartheta$ is constant along the connected
components of $\Pc$ defined by $\Dc$.
\end{prop}

\pf Let $A \in \B$. Then
\begin{equation}\label{squared}
\mathcal S^2(A) =
\dfrac{\cx(A)\cx(A^{-1})}{\cx(A^h)\cx(A^v)}A = \dfrac{\vartheta(bl(A))\vartheta(tr(A))}{\vartheta(br(A))\vartheta(tl(A))}A. \end{equation}
It follows that $\mathcal S^2({}_D\uno) = \frac{\vartheta(e(D))}{\vartheta(s(D))} {}_D\uno$, $D \in \Dc$. This proves the proposition. \epf

\bigbreak
\begin{obs}
Let $D \in \Dc$ and $E \in \Ec$. Then, as in Lemma \ref{regular},

\begin{equation}\label{antipoda-unos}
\Ss(\uno_E) = {}_{E^{-1}}\uno, \qquad \Ss({}_D\uno) =
\frac{\vartheta(e(D))}{\vartheta(s(D))}\, \uno_{D^{-1}}.
\end{equation}
\end{obs}

\bigbreak
\begin{prop}\label{pivotal}\footnote{We thank the referee for suggesting the formula for the pivotal element included here.}
There is a pivotal  group-like element $G\in \ku^{(\vartheta)} \T$
implementing $\Ss^2$ by conjugation-- see \cite{bnsz, nik-ss}.
Explicitly,
$$G = \sum_{x\in \Hc} \frac{\vartheta(l(x))}{\vartheta(r(x))}\,
 \id x,$$
where $\Hc(P, Q)$ is the set of all horizontal arrows going from $P$ to $Q$.
\end{prop}

Notice however that neither the elements $\uno_{\Theta(P)}$, $P\in \Pc$, nor the boxes $\id x$, $x\in \Hc$, are central in $\ku^{(\vartheta)} \T$.

\pf Note that $G = \Ss (w) w^{-1}$, where
$w = \sum_{P\in \Pc} \dfrac{1}{\vartheta(P)}\, \uno_{\Theta(P)}$.
Indeed,
$$
w^{-1} = \sum_{P\in \Pc} \vartheta(P)\, \uno_{\Theta(P)}
\qquad \text{and} \qquad
\Ss(w) = \sum_{Q\in \Pc} \frac{1}{\vartheta(Q)}\, {}_{\Theta(Q)}\uno,$$
the first equality by Lemma \ref{gpd-st}
and the second by \eqref{antipoda-unos}. Hence

\begin{align*}
G &= \sum_{P, Q\in \Pc} \frac{\vartheta(P)}{\vartheta(Q)}\,
\uno_{\Theta(P)}.{}_{\Theta(Q)}\uno \\
&= \sum_{P, Q\in \Pc} \frac{\vartheta(P)}{\vartheta(Q)}\,
\sum_{x\in \Hc(P, Q)} \id x\\
&= \sum_{x\in \Hc} \frac{\vartheta(l(x))}{\vartheta(r(x))}\,
 \id x.
\end{align*}
A straightforward computation shows that
$\Ss^2(A) = G^{-1}.A.G$, $A\in \B$. Finally,

\begin{align*}
\Delta(G) &= \Delta(1)\sum_{x\in \Hc} \quad \sum_{z,w\in \Hc: zw = x}
\frac{\vartheta(l(x))}{\vartheta(r(x))}\,
 \id z\otimes \id w
\\&= \Delta(1) \sum_{z,w\in \Hc}
\frac{\vartheta(l(z))}{\vartheta(r(z))}\,
\frac{\vartheta(l(w))}{\vartheta(r(w))}\, \id z\otimes \id w
\\ &= \Delta(1)(G \otimes G).
\end{align*}
Here the first equality is by  \eqref{delta-x}, and the second uses
that $l(z) = l(x)$, $r(w) = r(x)$ and $r(z) = l(w)$ whenever $x = zw$.
\epf

\bigbreak
\subsection{$C^*$-structure}

\

In this subsection $\ku = \mathbb C$.
Let $\vartheta: \Pc \to \mathbb R^+$ be a function satisfying
condition \eqref{multiplicativa-cx}, and let $\mathbb
C^{(\vartheta)}\T$ be the corresponding weak Hopf algebra. We
describe a $C^*$-structure on $\mathbb C^{(\vartheta)}\T$. See
\cite{bnsz}.

\medbreak Consider the map $\lambda: \B \to \mathbb R^+$, given by
$$\lambda(A) = \frac{\daleth(A^h)}{\daleth(A^{-1})} =
\frac{\vartheta(br(A))}{\vartheta(rt(A))}.$$

\begin{lema}\label{char} (i) $\lambda$ defines a character on the vertical composition groupoid
$\B \rightrightarrows \Hc$.

(ii) $\lambda$ is invariant with respect to horizontal
composition, that is, $\lambda(UV) = \lambda(V)$, for all $U, V
\in \B$ such that $U\vert V$.\end{lema}

\pf (i) Let $A, B \in \B$ such that $\frac{A}{B}$. We have $br(A)
= rt(B)$. Therefore,
$$\lambda(A) \lambda(B) =
\frac{\daleth(A^h)}{\daleth(A^{-1})}\frac{\daleth(B^h)}{\daleth(B^{-1})}
=
\frac{\vartheta(br(A))}{\vartheta(rt(A))}\frac{\vartheta(br(B))}{\vartheta(rt(B))}
= \frac{\vartheta(br(B))}{\vartheta(rt(A))}.$$ Since
$br\left(\begin{tiny}\begin{matrix} A \vspace{-4pt}\\B
\end{matrix}\end{tiny}\right) = br(B)$ and
$rt\left(\begin{tiny}\begin{matrix} A \vspace{-4pt}\\B
\end{matrix}\end{tiny}\right) = rt(A)$, the equality
$\lambda\left(\begin{tiny}\begin{matrix} A \vspace{-4pt}\\B
\end{matrix}\end{tiny}\right) = \lambda(A)\lambda(B)$ follows.
This proves part (i).

(ii) Let $U, V \in \B$ such that $U\vert V$. We have $br(A) =
rt(B)$. Hence, $\lambda(UV) =
\frac{\vartheta(br(UV))}{\vartheta(rt(UV))} =
\frac{\vartheta(br(V))}{\vartheta(rt(V))} = \lambda(V)$,  proving
part (ii). \epf

\begin{prop}\label{c*} $\mathbb C^{(\vartheta)}\T$ is a $C^*$-quantum groupoid with respect to the involution
$A^* : = \lambda(A)A^v$. \end{prop}

\pf We have $\lambda(A^v) = \lambda(A)^{-1}$. This implies that
$*$ is an involution, since $\lambda$ takes real values. Condition
$(A.B)^* = B^*.A^*$ follows from Lemma \ref{char}(i) and the fact
that \begin{tiny}$\left(\begin{matrix} A \vspace{-4pt}\\B
\end{matrix}\right)^v = \left(\begin{matrix} B^v
\vspace{-4pt}\\A^v \end{matrix}\right)$\end{tiny} for composable
boxes $A$, $B$. Also, for all $A \in \B$,
\begin{align*}\Delta(A^*) & = \lambda(A) \sum_{XY = A^v}
\frac{1}{\daleth(Y)} X \otimes Y = \lambda(A) \sum_{UV = A}
\frac{1}{\daleth(V^v)} U^v \otimes V^v \\ & = \sum_{UV = A}
\frac{\lambda(V)}{\daleth(V^v)} U^v \otimes V^v = \sum_{UV = A}
\frac{1}{\lambda(U) \daleth(V^v)} U^* \otimes V^*,
\end{align*} by Lemma \ref{char}(ii).
Since $rt(U) = lt(V)$ and $br(U) = bl(V)$, we have
$$\lambda(U)\daleth(V^v) = \frac{\vartheta(br(U)) \,
\vartheta(tl(V))}{\vartheta(rt(U))} =  \vartheta(br(U)) =
\daleth(V).$$ Thus $\Delta(A^*) = \Delta(A)^{* \otimes *}$.

\medbreak Let $\varphi: \B \to \mathbb R^+$ be given by
$\varphi(A) = \begin{cases}1, \quad \text{if} \; A \in \Hc,\\ 0,
\quad \text{otherwise},  \end{cases}$ $A \in \B$. Let also $(\quad
\vert \quad): \mathbb C^{(\vartheta)}\T \times \mathbb
C^{(\vartheta)}\T \to \mathbb C$ be the unique inner product such
that $(A\vert B) = \varphi(A^*.B)$, $A, B \in \B$. This makes
$\mathbb C^{(\vartheta)}\T$ into a $C^*$-quantum groupoid. \epf

\bigbreak
\subsection{Duality}

\

Let $\vartheta: \Pc \to \ku^{\times}$ be a generalized corner function satisfying \eqref{multiplicativa-cx}, and let $\ku^{(\vartheta)}\T$ be the corresponding weak Hopf algebra. We
describe in this subsection the dual weak Hopf algebra $(\ku^{(\vartheta)}\T)^*$.

\medbreak Consider the map $\mu: \B \to \ku^{\times}$, defined by
$$\mu(A) = \frac{\daleth(A^v)}{\daleth(A^{-1})} =
\frac{\vartheta(lt(A))}{\vartheta(rt(A))}.$$ The following properties are analogous to those in Lemma \ref{char}.

\begin{lema} (i) $\mu$ defines a character on the horizontal composition grou\-poid $\B \rightrightarrows \Vc$.

(ii) $\mu$ is invariant with respect to vertical
composition. \qed \end{lema}

Recall that the transpose double groupoid  $\T^t = \begin{matrix} \B &\rightrightarrows &\Vc
\\\downdownarrows &&\downdownarrows \\ \Hc &\rightrightarrows &\Pc \end{matrix}$ is obtained from $\T$ interchanging 'horizontal' and 'vertical' throughout. For $A \in \B$, let $A^t$ denote the same box in the transpose double groupoid; so that $$A^t = \begin{matrix} \quad l(A) \quad \\ t(A) \,\, \boxe \,\, b(A) \\ \quad r(A)\quad
\end{matrix}.$$

\begin{prop} $(\ku^{(\vartheta)}\T)^* \simeq \ku^{(\vartheta)}\T^{t}$ as weak Hopf algebras. \end{prop}

\pf Recall the basis $\{\underline{A}\}_{A \in \B}$ of $\ku\T$
 introduced in Remark \ref{nopref}.
Then the map $(\, , \, ): \ku^{(\vartheta)}\T \otimes \ku^{(\vartheta)}\T^{t} \to \ku$ defined by $(\underline{A}, B) : = \mu(A)\delta_{A, B^t}$ is a non-degenerate weak Hopf algebra pairing. \epf

\bigbreak
\subsection{Fusion categories}

\

We next establish necessary and sufficient conditions for
$\Rep \ku^{(\vartheta)} \T$ to be a fusion category.
For this, we study the $\ku^{(\vartheta)} \T$-module structure
of the target subalgebra
$\ku^{(\vartheta)} \T_t$. We begin by the following general result.

\bigbreak
 Let $\G$ be a groupoid with base $\Pc$ and let $\sim$ be the corresponding equivalence relation on $\Pc$: for $P, Q\in \Pc$,
$P\sim Q$ iff there exists $g\in \G$ such that $s(g) = P$, $e(g) = Q$.

\bigbreak Assume that $\G$ acts on a
fiber bundle $\gamma: \Ee \to \Pc$ via $\rightharpoondown:
\G \Times{e}{\gamma} \Ee \to \Ee$. If $P\in \Pc$, then we denote by $\Ee_P$ the fiber $\gamma^{-1}(P)$.
Clearly, $\gamma(\Ee)$ is stable under the equivalence relation $\sim$;
for, the action $g\rightharpoondown \quad: \Ee_{e(g)} \to \Ee_{s(g)}$
is a bijection.
We say that the action
is \emph{transitive} if $\gamma(\Ee)$ is just one class for the relation $\sim$.

\bigbreak
Let $\lambda: \G \to \ku^{\times}$ be a character, \emph{i. e.}
$\lambda(gh) = \lambda(g)\lambda(h)$ when $e(g) = s(h)$.
Let $\ku_{\lambda}\Ee$ be the $\ku\G$-module with $\ku$-basis $(x_E)$, $E\in \Ee$, and action
$$
g. x_E = \begin{cases} \lambda(g) x_{g \rightharpoondown E} \quad
&\text{if } e(g) = \gamma(E), \\ 0 \quad &\text{if not.}
\end{cases}
$$

\bigbreak
\begin{lema}\label{simplegeneral} The following are equivalent.

(i) The $\ku\G$-module $\ku_{\lambda}\Ee$ is simple.

(ii) The action is transitive and $\# \Ee_P = 1 $ for any
$P\in \gamma(\Ee)$.
\end{lema}

\pf (i) $\implies$ (ii). If $P\in \gamma(\Ee)$, we set
$$
x_P = \sum_{E\in \Ee_P} x_E.
$$

If $g\in \G$ has $s(g) = P$, $e(g) = Q$, then $g.x_Q = \lambda(g) x_P$. Let $\Q \subset \Pc$ be an equivalence class of $\sim$ and let $M$ be the $\ku$-span of the elements $x_P$, $P\in \Q$. Then
$$\oplus_{P\in \Q} \ku x_P = M = \ku_{\lambda}\Ee = \oplus_{P\in \Pc} \ku\Ee_P,$$
the second equality since $\ku_{\lambda}\Ee $
is simple. This implies that $\Q = \Pc$, \emph{i. e.} that the action
is transitive; and that $\dim \ku \Ee_P = 1$, \emph{i. e.} that  $\# \Ee_P = 1$, for any $P\in \gamma(Q)$.

\medbreak
(ii) $\implies$ (i). Let $M$ be a non-zero $\ku\G$-submodule of
$\ku_{\lambda}\Ee$ and let $0\neq m = \sum_{E\in \Pc} m_E x_E \in M$, where $m_E\in \ku$. Fix $E$ such that $m_E\neq 0$; then $x_E = m_E^{-1} \id \gamma(E) . m \in M$, hence $\ku\Ee_{\gamma(E)} \subset M$ by the assumption ``$\# \Ee_P = 1 $ for any $P\in \gamma(\Ee)$".
Now the transitivity assumption implies that
$\ku\Ee \subset M$. Hence $\ku_{\lambda}\Ee$ is simple.
\epf

\bigbreak Recall that a semisimple finite tensor category is called a
fusion category exactly when the unit object is simple \cite{ENO}.
A semisimple weak Hopf algebra is called \emph{connected} if its representation category is fusion \cite[Section 4]{ENO}.

\begin{prop}\label{prop:fusion}
The tensor category $\Rep \ku^{(\vartheta)} \T$ is a
fusion category (or $\ku^{(\vartheta)} \T$ is connected)
if and only if the following hold:
\begin{enumerate}
\item[(a)] $\Vc \rightrightarrows \Pc$ is connected, \emph{cf.}
Remark \ref{vertconn}.

\item[(b)]For any $x\in \Hc$, there exists at most one $E\in \Ec$
such that $b(E) = x$.
\end{enumerate}
\end{prop}

\pf Let $A\in \B$ and $E \in \Ec$. Then the action of $A$ on
$\ku \T_t$ is given in terms of the action $\curvearrowright$ in
\eqref{curve}. Explicitly,
\begin{multline*}
A \cdot\uno_E = \epsilon_t(A .\uno_E)= \epsilon_t(E
\rightharpoondown A) \\ = \begin{cases} \dfrac{\cx((E
\rightharpoondown A)^{-1})}{\cx((E
\rightharpoondown A)^h)}\uno_{\psi(E \rightharpoondown A)} = \dfrac{\cx(A^{-1})}{\cx(A^h)} \uno_{A \curvearrowright E},
&\text{if } b(E \rightharpoondown A) \in\Pc,
\\ 0 \quad &\text{if not.}
\end{cases}
\end{multline*}
Note that $b(E \rightharpoondown A)= b(E)b(A)$, thus
$b(E \rightharpoondown A) \in\Pc$ iff $b(E) = b(A)^{-1}$.

\medbreak
We can then apply Proposition \ref{simplegeneral} with $\gamma(E)
= b(E)^{-1}$, $E\in \Ec$.
The action $\curvearrowright$ is transitive iff (a) holds, by
Lemma \ref{simvsiisimcurve}; and condition (b) is equivalent to
``$\# \Ec_x = 1 $ for any $x\in \gamma(\Ec)$". \epf

In particular, if $\Rep\ku^{(\vartheta)} \T$ is fusion then necessarily $\Pc$ is vertically connected, because of condition (a) applied to the boxes
$\Theta(P)\in \Ec$.

\bigbreak
\subsection{Frobenius-Perron dimensions of simple objects}

\

Let $\T$ be a finite double groupoid and $\vartheta$ a generalized corner function satisfying \eqref{multiplicativa-cx}. Let
$\Rep\ku^{(\vartheta)} \T$ be the corresponding weak Hopf algebra.

\bigbreak
Let $\B \rightrightarrows \Hc$ be the vertical composition groupoid. Let $\R$ be the equivalence relation on $\Hc$ defined by $\B$ and let
$X^2$ be the coarse groupoid on $X\in \R$.
We fix $x\in X$, for  $X\in \R$, and denote by $\B(x) = \B(x,x)$
the group of loops in $x$. If $y\in \Hc$ belongs to the class $X$,
then we denote $\overline{y} = x$.
The groupoid structure of $\B$ is determined by

\begin{equation}\label{estructuragpde}
\B \simeq \coprod_{X\in \R} X^2 \times \B(x).
\end{equation}

If $F$ is a finite set, we denote by $\ku F$ the vector space
with base $F$ and by
$M_F (\ku)$ the `matrix algebra' $\End (\ku F)$,
with matrix idempotents $E_{st}$, $s,t\in F$.
Then \eqref{estructuragpde} gives in turn an isomorphism of algebras
$$
\ku \T \simeq \prod_{X\in \R} M_X (\ku) \otimes \ku\B(x)
$$
that maps the identity arrows to primitive idempotents:
$$
\id y \mapsto E_{yy} \otimes e \in M_X (\ku) \otimes \ku\B(x),
$$
$y\in X$, $X\in \R$.
Clearly any simple $\ku^{(\vartheta)} \T$-module is of the form
$U = \ku X\otimes V$, where $X\in \R$ and $V$ is a simple
$\B(x)$-module.

\bigbreak
\emph{Suppose for the rest of this subsection
that $\Rep \ku^{(\vartheta)} \T$ is a fusion category
(or that $\ku^{(\vartheta)} \T$ is connected),
see Proposition \ref{prop:fusion}; and that $\ku = \mathbb C$}.
Recall the pivotal group-like element $G$ computed in
Proposition \ref{pivotal}.
Then $\Rep \ku^{(\vartheta)} \T$ is a pivotal fusion category
\cite[Definition 2.7]{ENO} and the quantum dimension
(corresponding to this pivotal structure) of the
simple module $U = \ku X\otimes V$ is
$$\qdim U = \dfrac{\tr_U (G)}{\dim \ku^{(\vartheta)} \T_t},$$  thus
\begin{equation}\label{formula:qdim}
\qdim U = \dfrac{1}{\#\Ec} \sum_{y\in X}
\dfrac{\vartheta(l(y))}{\vartheta(r(y))} \tr_U (\id y) =
\dfrac{\dim V}{\#\Ec} \sum_{y\in X}
\dfrac{\vartheta(l(y))}{\vartheta(r(y))}.
\end{equation}

Therefore, the \emph{quantum dimension} $\qdim \ku^{(\vartheta)} \T$ is
given by

\begin{align*}\qdim \ku^{(\vartheta)} \T  & = \sum_{X\in \R, V\in \widehat{\B(x)}} \Big(\dfrac{\dim V}{\#\Ec}\Big)^2
\Big\vert  \sum_{y\in X}
\dfrac{\vartheta(l(y))}{\vartheta(r(y))}\Big\vert^2 \\
& = \dfrac{1}{(\#\Ec)^2}\sum_{X\in \R}
\Big\vert  \sum_{y\in X}
\dfrac{\vartheta(l(y))}{\vartheta(r(y))}\Big\vert^2
|\B(x)|,\end{align*}
\emph{cf.} \cite[Definition 2.2 and Proposition 2.9]{ENO}.

\begin{lema}\label{lema:pfrob} If
$\sum_{y\in X}
\dfrac{\vartheta(l(y))}{\vartheta(r(y))} > 0 $
for any class $X\in \R$, then the
Frobe\-nius-Perron dimensions of the simple
$\ku^{(\vartheta)} \T$-modules agree with their quantum dimensions,
hence they are given by \eqref{formula:qdim}.
\end{lema}

\pf This follows because $\fpd$ is the unique ring homomorphism
from the Grothendieck ring to the complex numbers taking positive real values on the irreducible modules. \epf

When $\vartheta$ takes values in $\mathbb Q_{>0}$
(for example, when
$\vartheta = \dfrac 1{\theta}$ comes from the corner function),
Lemma \ref{lema:pfrob} says that the
Frobenius-Perron dimension of a simple
$\ku^{(\vartheta)} \T$-module $U$ is a rational number;
since by \cite{ENO}
$\fpd U$ is an algebraic integer, it is an integer. That is,
the Frobenius-Perron dimension of the
simple $\ku \T$-modules are integers:
$$
\dfrac{\dim V}{\#\Ec} \sum_{y\in X}
\dfrac{\theta(r(y))}{\theta(l(y))} \in \mathbb N,
$$
for any class $X\in \R$, for any irreducible representation $V$
of $\B(x)$. In particular $\Rep \ku \T$
is tensor equivalent to the representation category of a finite-dimensional semisimple quasi-Hopf algebra \cite[Theorem 8.33]{ENO}.

\bigbreak
In general, a connected weak Hopf algebra $H$ is called
\emph{pseudo-unitary}  if $\qdim H$ coincides with the
Frobenius-Perron dimension of $H$.

\begin{prop}\label{prop:psunit} Assume that
\begin{equation}\label{formula:psunit}
\dfrac{\cx(A)\cx(A^{-1})}{\cx(A^h)\cx(A^v)} = \dfrac{\vartheta(lb(A))\vartheta(rt(A))}{\vartheta(rb(A))\vartheta(lt(A))} > 0
\end{equation}
for any box $A\in \B$. Then
$\ku^{(\vartheta)} \T$ is a pseudo-unitary weak Hopf algebra and
for any irreducible $\ku^{(\vartheta)} \T$-module $U = \ku X \otimes V$ as above,
$$
\fpd U = \dfrac{\dim V}{\#\Ec} \Big\vert  \sum_{y\in X}
\dfrac{\vartheta(l(y))}{\vartheta(r(y))}\Big\vert.$$
\end{prop}

\pf By \eqref{squared}, condition \eqref{formula:psunit}
says that all eigenvalues of $\mathcal S^2$ are strictly positive
numbers; the first claim follows from the criterion given in
\cite[Corollary 5.2.5]{nik-ss}.  The second claim follows from
\cite[Proposition 8.21]{ENO}.  \epf

\section{Examples}\label{ejemplos}

\subsection{Matched pairs}\label{ap-tak}\textbf{\cite{AN}.}

\

Given a groupoid $\G$  one can consider the double groupoid whose
boxes are commuting squares in $\G$. More generally, if $\Hc$ and
$\Vc$ are wide subgroupoids  of $\G$, there is a double groupoid
of commuting squares in $\G$ whose horizontal arrows belong to
$\Hc$ and whose vertical arrows belong to $\Vc$. A special case of
this remark is given by the following construction.

Let $\fiz: \Hc {\,}_r \times_t \Vc \to \Hc$, $\fde: \Hc {\,}_r
\times_t \Vc \to \Vc$ be a matched pair of finite groupoids, on
the same basis $\Pc$. Here, $l, t$ (respectively, $r, b$) denote
the source and target maps of $\Hc$ and $\Vc$. This is equivalent
to giving an exact factorization $\G = \Vc {}_b\times_{l}\Hc$.

Let $\B: = \Vc \Times{t}{r} \Hc$. We have a double groupoid
$$\begin{matrix} \Vc \Times{t}{r} \Hc &\rightrightarrows & \Hc
\\\downdownarrows &&\downdownarrows \\ \Vc &\rightrightarrows &\Pc, \end{matrix}$$
defined as follows:  $\Vc \Times{t}{r} \Hc \rightrightarrows  \Hc$
and $\Vc \Times{t}{r} \Hc \rightrightarrows \Vc$ are,
respectively, the transformation groupoids corresponding to the
actions $\fiz$ and $\fde$. These examples exhaust the class of
vacant double groupoids, as shown by Mackenzie \cite{mk1}. The
corresponding weak Hopf algebras have been introduced in \cite{AN}.

The core groupoid of a vacant double groupoid is the discrete
groupoid on $\Pc$. Therefore the source and target subalgebras of
the associated quantum groupoid are commutative, {\it i.e.}, this
quantum groupoid is a \emph{face algebra} in the sense of Hayashi \cite{H}.

\medbreak For Hopf algebras coming from matched pairs of finite groups, the pictorial description using boxes has been used in \cite{maj-lib, tak}.

\bigbreak
\subsection{Bimodules over a separable algebra}\label{bimodulos}

\

Let $\G$ be a finite groupoid with source and target maps $s, e:
\G \rightrightarrows \Pc$. There is a double groupoid
$$\T = \T(\G) = \begin{matrix} \G \times \G &\rightrightarrows &\Pc \times
\Pc \\\downdownarrows &&\downdownarrows \\ \G &\rightrightarrows
&\Pc \end{matrix},$$ canonically associated to $\G$, where
\begin{itemize}
    \item $\Pc \times \Pc \rightrightarrows \Pc$ is the coarse groupoid on
$\Pc$,
    \item $\G \times \G \rightrightarrows \G$ is the coarse groupoid on
$\G$,
    \item $\G \times \G \rightrightarrows \Pc \times \Pc$ is the product
    groupoid of $\G$ with itself.
\end{itemize}
See \cite[Example 1.3]{BM}. A box in this double groupoid is a
pair $(g, h)$ of arrows in $\G$ and it can be depicted as $$(g, h)
=
\begin{matrix} \quad \quad
\\  g\,\, \boxe  \,\, h \\ \quad  \quad
\end{matrix} =  \begin{matrix} s(g) \quad
s(h) \\  g \,\, \boxed{(g, h)} \,\, h \\ e(g) \quad e(h)
\end{matrix}.$$ Horizontal and vertical compositions are  defined,
respectively, by \begin{equation*} \begin{matrix} \quad  \quad \\
g\,\, \boxe \,\, h \\ \quad  \quad
\end{matrix} \quad  \begin{matrix} \quad  \quad \\  h\,\, \boxe  \,\, t \\ \quad  \quad
\end{matrix}  = \begin{matrix} \quad  \quad \\  g \,\, \boxe  \,\, t \\ \quad  \quad
\end{matrix}, \qquad
\begin{matrix} \quad  \quad \\  g\,\, \boxe  \,\, h \\ \quad  \quad \\  g' \,\, \boxe  \,\, h' \\ \quad  \quad \end{matrix}  =
\begin{matrix}  \quad  \quad \\  gg'\,\, \boxe  \,\, hh' \\ \quad  \quad
\end{matrix},\end{equation*} for all pairs of composable arrows $g, g'$ and $h, h'$, and for all arrow $t$.
The horizontal
and vertical identities are given by
$$\iddv g = \begin{matrix} \quad  \quad \\  g\,\, \boxe  \,\, g \\ \quad  \quad
\end{matrix}, \qquad \iddv (p, q) = \begin{matrix} \quad  \quad \\  \id_p \,\, \boxe  \,\, \id_q \\ \quad  \quad
\end{matrix}, $$
for all $g \in \G$, $p, q \in \Pc$.

We summarize the relevant facts on the structure of this double
groupoid.

\begin{lema}\label{bs-bimod} (a). The core groupoid $\Dc$ is isomorphic to $\G$.

(b). Let $g, h \in \G$ and let $A = \begin{matrix}
\quad\quad\vspace{-4pt} \\  g\,\, \boxe  \,\,h   \\
\quad\quad\vspace{-4pt} \end{matrix} \in \B$. Then

(i) $\urcorner(A) = \# \{ u \in \G: \; s(u) = s(g) \}$.

(ii) $\llcorner(A) = \# \{ u \in \G: \; e(u) = e(h) \}$.

(iii) $\ulcorner(A) = \# \{ u \in \G: \; s(u) = s(h) \}$.

(iv) $\lrcorner(A) = \# \{ u \in \G: \; e(u) = e(g) \}$.\qed \end{lema}

By Lemma \ref{bs-bimod}, the source subalgebra of the
corresponding weak Hopf algebra $\ku \T(\G)$ is isomorphic to $\ku
\G$. Also, since $\urcorner(A) = \lrcorner(A)$, and $\llcorner(A)
= \ulcorner(A)$, for all $A \in \B$, $\mathcal S^2 = \id$ in $\ku
\T(\G)$, \emph{cf.} Lemma \ref{sq-antp}.

\begin{obs} Observe that in the double groupoid $\T(\G) \coprod {}^t\T(\G)$, where ${}^t\T(\G)$ is the transpose double groupoid, the  corner functions give four pairwise \emph{distinct} maps $\B \to \mathbb N$. \end{obs}

\begin{prop}\label{d(g)}The category $\Rep \ku \T(\G)$ is equivalent to the
category ${}_{\ku\G}\mathcal M_{\ku\G}$  of finite dimensional
$\ku \G$-bimodules. \end{prop}

\pf The algebra $\ku\G$ is a separable algebra, with separability
idempotent $$e = \sum_{g \in \G} \dfrac{1}{d(s(g))} \; g \otimes
g^{-1} \, \in \, \ku\G \otimes (\ku\G)^{\op},$$ where $d(s(g))$ is
the number of arrows of $\G$ having source $s(g)$. Hence, for all
$\ku\G$-bimodules $V$ and $W$, the inclusion $e.(V\otimes_{\ku}W)
\subseteq V\otimes_{\ku}W$ induces a natural isomorphism
$$V\otimes_{\ku\G}W \simeq e.(V\otimes_{\ku}W).$$

We show in what follows that this coincides with the tensor
product of $\ku\T(\G)$-modules. As an algebra, $\ku \T(\G)$ is
isomorphic to $\ku (\G \times \G)$; this is in turn isomorphic to
$\ku\G \otimes (\ku\G)^{\op}$, the isomorphism given by inverting
the second factor. This gives a natural equivalence between $\ku
\T(\G)$-modules and $\ku\G$-bimodules. Tensor product of two $\ku
\T(\G)$-modules $V$ and $W$ is defined by $V \otimes W =
\Delta(1).(V \otimes_{\ku} W)$. Explicitly, we have

\begin{align*}\Delta(1) & = \sum_{p, q \in \Pc, h \in
\G}\dfrac{1}{d(s(h))} \quad
\begin{matrix} \quad  \quad \\  \id_p\,\, \boxe  \,\, h \\ \quad
\quad \end{matrix} \otimes \begin{matrix} \quad  \quad \\  h\,\,
\boxe  \,\, \id_q \\ \quad  \quad \end{matrix} \\ & =
\sum_{h \in \G}\dfrac{1}{d(s(h))}
\quad \left(\sum_{p\in \Pc}
\begin{matrix} \quad  \quad \\  \id_p\,\, \boxe  \,\, h \\ \quad
\quad \end{matrix}\right) \otimes \left(\sum_{q \in \Pc}\begin{matrix} \quad  \quad \\  h\,\,
\boxe  \,\, \id_q \\ \quad  \quad \end{matrix}\right),\end{align*} which
corresponds to the element $e \in \ku\G \otimes (\ku\G)^{\op}$. Therefore $V\otimes W$ corresponds to $V\otimes_{\ku\G}W$ as linear spaces. Finally, the action of a box $\begin{matrix} \quad  \quad \vspace{-6pt}\\  g\,\, \boxe  \,\,h  \\
\quad  \quad \vspace{-6pt}\end{matrix}$ on this tensor product coincides with the action of the corresponding element in $\ku\G \otimes (\ku\G)^{\op}$. This proves the lemma.
\epf

Observe that any  algebra $R$ which is separable over $\ku$ is
isomorphic as an algebra, to the groupoid algebra of a (not
canonical) finite groupoid. The following proposition is a
consequence of Lemma \ref{d(g)}.

\begin{prop} Let $R$ be separable algebra over $\ku$.
There exists a finite groupoid $\G$ such that
${}_R\mathcal M_R$ is tensor equivalent to $\Rep \ku \T(\G)$. \qed
\end{prop}

\bigbreak
\subsection{Group theoretical fusion categories}

\

In this section we consider the group theoretical categories introduced by Ostrik in \cite{ostrik}. Let $\Vect : = \Vect_{\ku}$ denote the category of vector spaces over $\ku$.

Let $G$ be a finite group, and let $F \subseteq G$ be a subgroup.
There is a fusion category $\C(G, \omega, F, \alpha)$, called a
\emph{group theoretical} category \cite[Definition 8.46]{ENO},
associated to  the following data:
\begin{itemize}\item a normalized 3-cocycle $\omega: G \times G \times G \to k^{\times}$;
\item a normalized 2-cochain $\alpha: F \times F \to k^{\times}$;
\end{itemize}
subject to the condition \begin{equation}\label{rest-to-F}\omega\vert_{F\times F \times F} = d\alpha. \end{equation}

By \eqref{rest-to-F}, the twisted group algebra $k_{\alpha}F$ is
an (associative unital) algebra in $\Vect^G_{\omega}$ of
$G$-graded vector spaces with associativity defined by $\omega$.
The category $\C(G, \omega, F, \alpha)$ is by definition the
$k$-linear monoidal category
${_{\,}}_{\ku_{\alpha}F}\hspace{-4pt}\left(\Vect^G_{\omega}\right)_{\ku_{\alpha}F}$ of
$k_{\alpha}F$-bimodules in $\Vect^G_{\omega}$: tensor product is
$\otimes_{k_{\alpha}F}$ and the unit object is $k_{\alpha}F$.  A
(quasi)-Hopf algebra $A$ is called group theoretical if the
category $\Rep A$ of its finite dimensional representations is
group theoretical. Every group theoretical category is the
representation category of a semisimple finite dimensional
quasi-Hopf algebra.

It is shown in \cite{fs-indic} that there is a 3-cocycle
$\widetilde \omega$, cohomologous  to $\omega$, such that the
categories $\C(G, \omega, F, \alpha)$ and $\C(G, \widetilde
\omega, F, 1)$ are equivalent. That is, up to tensor
equivalence, it is enough to consider the case
$\alpha = 1$.

\bigbreak
\subsection{The category $\Vect^G_{\omega}$}

\

In the next example we shall consider  the double groupoid $$\T_0
= \quad
\begin{matrix} G \times G\times G &\rightrightarrows & G \times G
\\\downdownarrows &&\downdownarrows \\ G \times G &\rightrightarrows &
G \end{matrix},$$ where
\begin{itemize}
\item the horizontal groupoid $G \times G\times G\rightrightarrows G \times G$
is the direct product of the  coarse groupoid on $G$ and the
discrete groupoid on $G$; that is, the  double groupoid
corresponding to the equivalence relation defined on $G \times G$
by $(x, y) \sim (x', y')$ if and only if $y = y'$;
    \item the vertical groupoid
     $G \times G\times G \rightrightarrows G \times G$ is the transformation
    groupoid associated to the regular action $.: G \times G \times G    \to G \times G$,
    $(x, y) . g = (xg, yg)$;
    \item the horizontal groupoid $G \times G \rightrightarrows G$
    is the coarse groupoid on $G$;
    \item the vertical groupoid $G \times G \rightrightarrows G$
    is the transformation groupoid corresponding to the right regular action of $G$ on
    itself.
\end{itemize}
Observe that $\T_0$ is a \emph{vacant} double groupoid, with boxes
determined by
$$(a, b, g): =  \begin{matrix}\quad a  \quad b \quad
\\  \,\, \boxe  \,\,  \\ \quad ag  \quad bg \quad \end{matrix},$$ for all $a, b, g \in
G$. Horizontal and vertical compositions are as follows:
\begin{equation*}\begin{matrix}\quad a  \quad b \quad
\\  \,\, \boxe  \,\,  \\ \quad ag  \quad bg \quad \end{matrix} \,
\begin{matrix}\quad b  \quad c \quad
\\  \,\, \boxe  \,\,  \\ \quad bg  \quad cg \quad \end{matrix} =
\begin{matrix}\quad a  \quad c \quad
\\  \,\, \boxe  \,\,  \\ \quad ag  \quad cg \quad \end{matrix}, \qquad  \begin{matrix}\quad a  \quad b \quad
\\  \,\, \boxe  \,\,  \\  \quad ag  \quad bg \quad  \\  \,\, \boxe  \,\,
\\ \quad agh  \quad bgh \quad \end{matrix} = \begin{matrix}\quad a  \quad b \quad
\\  \,\, \boxe  \,\,  \\ \quad agh  \quad bgh \quad \end{matrix}. \end{equation*}

\medbreak Let $\sigma: \B \times_{b, t} \B \to \ku^{\times}$ be
given by
\begin{equation}\sigma\left( \begin{matrix}\quad a  \quad b \quad
\\  \,\, \boxe  \,\,  \\ \quad ag  \quad bg \quad \end{matrix},
\begin{matrix}\quad ag  \quad bg \quad
\\  \,\, \boxe  \,\,  \\ \quad agh  \quad cgh \quad \end{matrix} \right) =
\dfrac{\omega(a, g, h)}{\omega(b, g, h)}. \end{equation} The
3-cocycle condition on $\omega$ implies that $\sigma$ is a
normalized vertical 2-cocycle for $\T_0$; see \cite{AN}. Moreover,
for all boxes $A$, $B$, $C$,
$D$ such that $\begin{tabular}{p{0,4cm}|p{0,4cm}} $A$ & $B$ \\
\hline $C$ & $D$
\end{tabular}$, we have $\sigma(AB, CD) = \sigma(A, C)\sigma(B,
D)$. Thus $\sigma$ and the trivial cocycle $\tau = 1$ are
compatible in the sense of \cite[Definition 3.6]{AN}, and
therefore there is an associated weak Hopf algebra
$\ku_{\sigma}\T_0$.

\begin{prop} $\Rep \ku_{\sigma}\T_0 \simeq
\Vect^G_{\omega}$ as  tensor categories. \end{prop}

In particular, the Drinfeld center $\mathcal Z(\Vect^G_{\omega})$ is
equivalent to the representation category  of the quantum double
$D(\ku_{\sigma}\T_0)$.

\begin{obs} Compare this example with the weak Hopf algebra $\mathcal A^{\omega}G$ described in \cite[Appendix]{bsz}. \end{obs}

\pf As an
algebra, $\ku_{\sigma}\T_0 \simeq \ku^{G\times G}\#_{\sigma}\ku G$
is the crossed product corresponding to the right regular action
of $G$ on $G\times G$ and the cocycle $\sigma(g, h) = \sum_{a, b
\in G}\sigma^a_b(g, h)e^a_b$, $g, h \in G$, where $e^a_b\in
\ku^{G\times G}$ are the canonical idempotents and $\sigma^a_b(g,
h): = \dfrac{\omega(a, g, h)}{\omega(b, g, h)}$. The
comultiplication in the canonical basis of $\ku^{G\times
G}\#_{\sigma}\ku G$ is determined by
\begin{equation}\label{com-l}\Delta(e^a_b \# g) = \sum_{c \in G} e^a_c \# g
\otimes e^c_b \# g,\end{equation} for all $a, b, g \in G$.

Consider the quasi-Hopf algebra structure on $H = \ku^G$ with
associator $\omega$. The representation category of $H$ is exactly
$\Vect^G_{\omega}$. By the results of Hausser and Nill \cite{HN},
there is a tensor equivalence $\Vect^G_{\omega} \to {}_H\mathcal
M^H_H$, where ${}_H\mathcal M^H_H$ is the category of quasi-Hopf
bimodules with tensor product $\otimes_H$. In particular, the
forgetful functor ${}_H\mathcal M^H_H \to {}_H\mathcal M_H$ is
monoidal, hence a fiber funtor. The defining relations for the
category ${}_H\mathcal M^H_H$ can be used to reconstruct a
multiplication and comultiplication in the vector space $L =
\ku^{(G\times G)} \otimes \ku G$, which make $L$ into a weak Hopf
algebra with basis $H$ (hence a face algebra, since $H$ is
commutative)  such that ${}_H\mathcal M^H_H \simeq \Rep L$.
Moreover, the formulas thus obtained give the crossed product
algebra structure $L = \ku^{G\times G}\#_{\sigma}\ku G$ with
comultiplication \eqref{com-l}. This establishes the proposition.
\epf

\bigbreak
\subsection{The categories ${}_F\Vect^G_F$}

\

\emph{We assume in this subsection that the cocycle $\omega$ is
trivial.} Consider the double groupoid $$\T = \quad \begin{matrix}
F \times F\times G &\rightrightarrows & F \\\downdownarrows
&&\downdownarrows \\ G &\rightrightarrows & * \end{matrix},$$
where $*$ is a set with a single element, and
\begin{itemize}
    \item $F \times F\times G   \rightrightarrows F$ is the direct product
    of the  coarse groupoid on $F$ and group $G$,
    \item $F \times F\times G \rightrightarrows G$ is the transformation
    groupoid associated to the action $.: G \times F \times F    \to G$, $g.(x, y) = x^{-1}gy$.
\end{itemize}
The box in $\T$ corresponding to a triple $(x, y, g)$, $x, y \in F$, $g \in G$ can be depicted as
$$ \begin{matrix} \quad g \quad
\\  x\,\, \boxe  \,\, y \\ \quad h \quad
\end{matrix},$$
where $h \in G$, $gy = xh$.
Horizontal and vertical compositions are given by
\begin{equation*} \begin{matrix} \quad g \quad
\\  x\,\, \boxe  \,\, y \\ \quad h \quad
\end{matrix} \,  \begin{matrix} \quad s \quad
\\  y\,\, \boxe  \,\, z \\ \quad t \quad
\end{matrix}  = \begin{matrix} \quad gs \quad
\\  x\,\, \boxe  \,\, z \\ \quad ht \quad
\end{matrix}, \quad \begin{matrix} \quad g \quad
\\  x\,\, \boxe  \,\, y \\ \quad h \quad \\  \quad h\quad \\ x'\,\, \boxe  \,\,
y' \\ \quad t \quad \end{matrix}  =
\begin{matrix} \quad g \quad
\\  xx'\,\, \boxe  \,\, yy' \\ \quad t, \quad
\end{matrix},\end{equation*} for all $x, y, z, x', y' \in F$, $g, h, s, t \in G$.
In this example the core groupoid $\Dc$ is isomorphic to $F$,
and the right corner function is given by the
formula \begin{equation}\urcorner\left( \begin{matrix} \quad g
\quad \\  x\,\, \boxe  \,\, y \\ \quad h \quad
\end{matrix} \right) = \vert F\vert. \end{equation}
$\T$ is the \emph{comma double groupoid} associated to the
inclusion $F \to G$ \cite[Example 1.8]{BM}. In the terminology
{\it loc. cit.} $\T$ coincides with the transformation double
groupoid $\T(F) \ltimes (\chi_G, \id)$, where $\T(F)$ is the
double groupoid in Subsection \ref{bimodulos}, and  $\chi_G$ is
the anchor map corresponding to $G$.

\begin{prop} $\ku \T \simeq
{}_F\Vect^G_F$ as  tensor categories. \end{prop}

\pf Since the vertical groupoid $\B \rightrightarrows \Hc$ is a
transformation groupoid, $\ku \T \simeq \ku^G \# \ku(F \times
F^{\op})$ as an algebra, where the right hand side expression is
the smash product corresponding to the action $\ku(F \times
F^{\op}) \otimes \ku^G \to \ku^G$, $\langle(x, y).f, g\rangle :=
\langle f, g.(x, y)\rangle$. Then there is an equivalence $\ku\T
\simeq {}_F\Vect^G_F$ of $\ku$-linear categories.
Under this identification, the comultiplication on $\ku \T$ is
given by
\begin{equation*}\Delta(e_g \# (x \otimes y)) = \sum_{st = g}
\left(e_s \# x \otimes e^{(1)}\right) \otimes \left(e_t \# e^{(2)}
\otimes y\right), \end{equation*} where $e_g \in \ku^G$ are the
canonical idempotents, $g \in G$, and $e = e^{(1)} \otimes e^{(2)}
\in \ku F \otimes (\ku F)^{\op}$ is the separability idempotent
given by $e = \sum_{x \in F} \dfrac{1}{\vert F\vert} x \otimes
x^{-1}$.

In particular, $\Delta(1) = \sum_{s, t} \left(e_s \# 1_G \otimes
e^{(1)}\right) \otimes \left(e_t \# e^{(2)} \otimes 1_G\right)$.
Therefore, as in the proof of Lemma \ref{d(g)}, we see that the
natural equivalence $\ku\T \simeq {}_F\Vect^G_F$ of $\ku$-linear
categories preserves tensor products. This finishes the proof of
the proposition.  \epf


\end{document}